\documentclass{amsart}

\newtheorem{theo}{Theorem}[section]

\newtheorem{theorem}[theo]{Theorem}

\theoremstyle{definition}

\theoremstyle{remark}
\newtheorem{remark}[theo]{Remark}
\newtheorem{example}[theo]{Example}

\usepackage{tikz}
\usetikzlibrary{cd}
\usepackage{here}
\usepackage{amssymb}
\usepackage{enumitem}
\usepackage{scalerel}
\usepackage{url}
\usepackage{cite}

\title{On a construction of stable maps from $3$-manifolds into surfaces}

\author{Gakuto Kato}
\address{Graduate School of Integrated Basic Sciences, Nihon University, 3-25-40 Sakurajosui, Setagaya-ku, Tokyo 156-8550, JAPAN}
\email{chga24002@g.nihon-u.ac.jp}

\date{\today}

\begin{document}

\keywords{stable map, $3$-manifold, braid}

\subjclass[2020]{57R45, 57M99, 57K10}

\begin{abstract}
For any link in the $3$-sphere, we provide a visual construction of a stable map $f$ from the $3$-sphere to the Euclidean plane such that $f$ has no cusp points, the set of definite fold points of $f$ is isotopic to the given link, and $f$ has only a certain type of singular fibers. As a corollary, we obtain a stable map $f_{0}$ from every closed orientable smooth $3$-manifold to the $2$-sphere such that $f_{0}$ has neither cusp points nor definite fold points, and has only a certain type of singular fibers.
\end{abstract}

\maketitle

%%%%%%%%%%%%%%%%%%%%%%%%%%%%%%%%%%%%%%%%%%%%%

\section{Introduction} 
\label{s1}

There have been many studies focusing on stable maps from a closed orientable smooth $3$-manifold to the Euclidean plane $\mathbb{R}^{2}$. See \cite{Hiratuka-Jorge-Saeki, Ishikawa-Koda, Kalmar-Stipsicz, K.-Levine-Port., Levin'65, Saeki1993, Saeki1994, Saeki'95, Saeki'96, Saeki2019, Saeki-Yamamoto2016, Saeki-Yamamoto2018, S.Naoki}, for example. The definition of stable maps will be introduced in Section~\ref{s2}. In this paper, all manifolds and maps are assumed to be differentiable of class $C^\infty$ unless otherwise indicated.

\begin{figure}[htbp]
        {\unitlength=1cm
        \begin{picture}(12.5,1.5)(0,0)
        \put(4,-0.5){\includegraphics[height=3cm,clip]{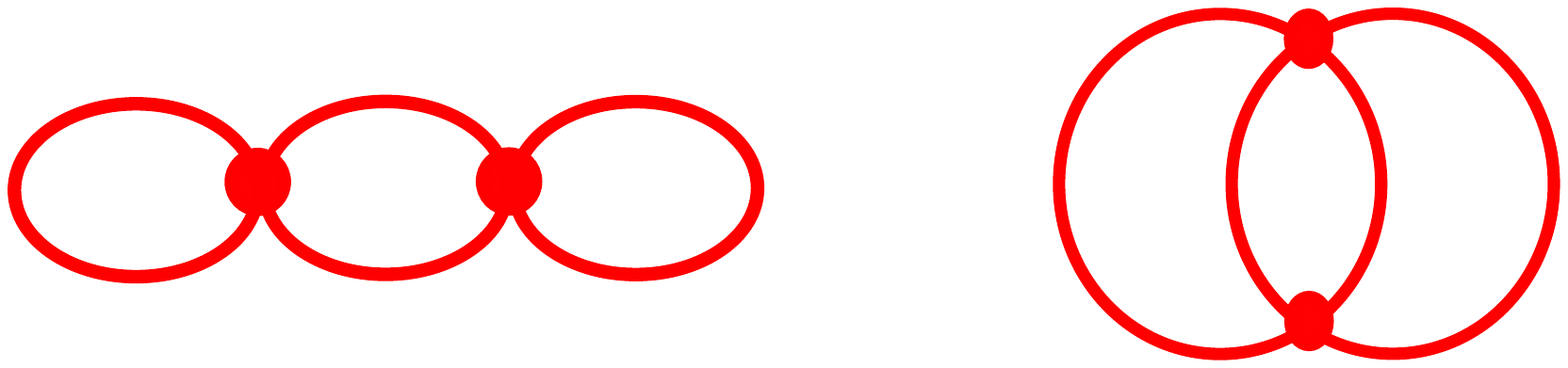}}
        \put(4.5,0.2){type $\mathrm{I\hspace{-1.2pt}I^{2}}$}
        \put(7,0.2){type $\mathrm{I\hspace{-1.2pt}I^{3}}$}
        \end{picture}}
        \caption{Types of singular fibers.}
        \label{types of sing. fibers}
    \end{figure}

Let $M$ be a closed orientable $3$-manifold and let $L$ be any given link in $M$. It is shown by Saeki \cite[Corollary 6.3]{Saeki'95, Saeki'96} that there exists a stable map $f : M \to \mathbb{R}^{2}$ such that the set of singular points $S(f)$ is equal to $L$ if and only if $[L]_{2}=0$ in $H_{1} (M;\mathbb{Z}_{2})$. For such a stable map $f$, it is known that the points in $S(f)$ are classified into definite fold points, indefinite fold points, and cusp points. Further, it is known that singular fibers of $f$ containing two indefinite fold points are classified into \textit{type $\mathrm{I\hspace{-1.2pt}I^{2}}$} and \textit{type $\mathrm{I\hspace{-1.2pt}I^{3}}$} (Figure~\ref{types of sing. fibers}). See \cite{Levine1985, Levin'65}, for example. Then, in \cite[Corollary 3.7]{Ishikawa-Koda}, Ishikawa and Koda showed that there exists a stable map $f \colon M \to \mathbb{R}^{2}$ without cusp points such that $f$ has no singular fibers of type $\mathrm{I\hspace{-1.2pt}I^{3}}$ and the set of definite fold points $S_{0}(f)$ contains the given link $L$.

% When $y \in \mathbb{R}^{2}$ is a regular value of $f$, we call $f^{-1}(y)$ a \textit{regular fiber}; otherwise, it is called a \textit{singular fiber} ( for detail, see \cite{Saeki04}).

    \begin{figure}[htbp]
        {\unitlength=1cm
        \begin{picture}(12.5,4)(0,0)
        \put(3,0){\includegraphics[height=5cm,clip]{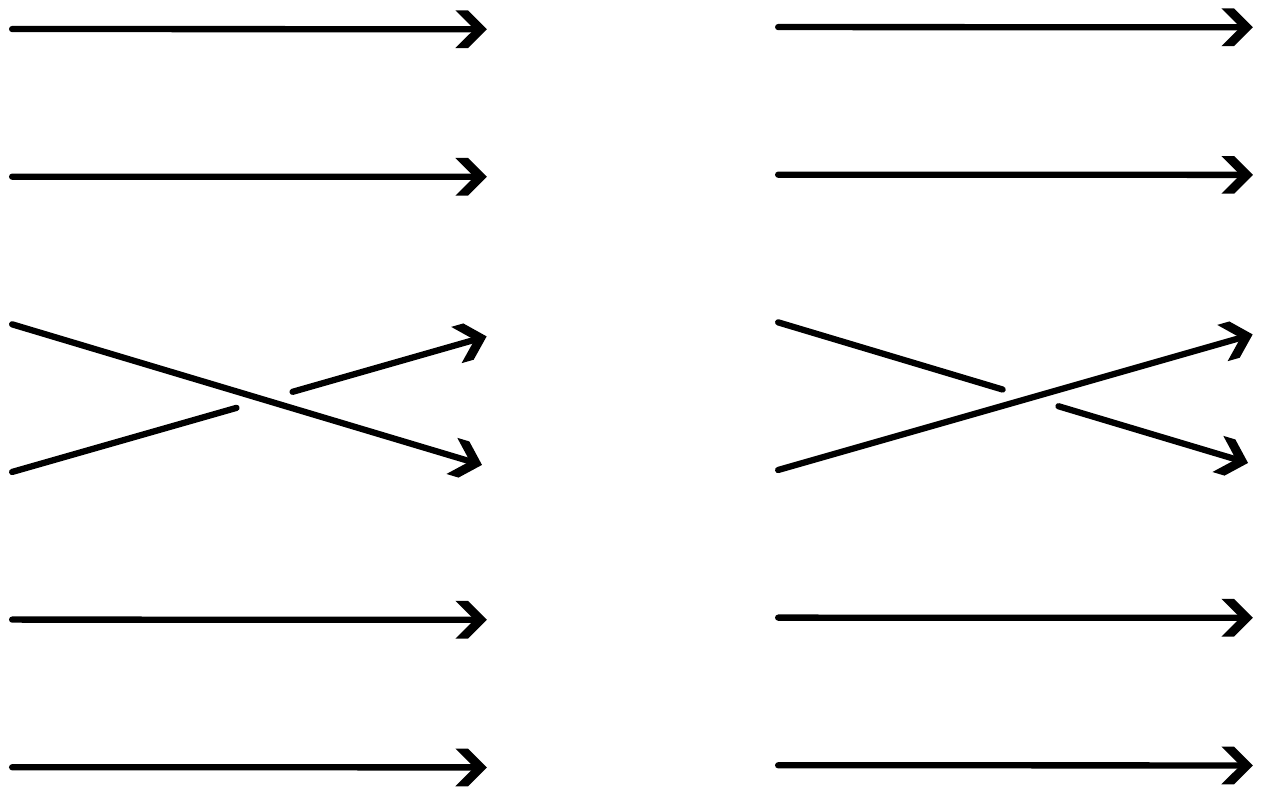}}
        \put(3.5,3.8){$1$}
        \put(3.5,3.2){$2$}
            \put(4.9,2.8){$\vdots$}
        \put(3.5,2.7){$i$}
        \put(3,2){$i+1$}
            \put(4.9,1.8){$\vdots$}
        \put(3,1.4){$n-1$}
        \put(3.5,0.9){$n$}
            \put(8,2.8){$\vdots$}
            \put(8,1.8){$\vdots$}
        
        \put(4.8,0.4){$\sigma_{i}$}
        \put(7.7,0.4){$\sigma_{i}^{-1}$}
        \end{picture}}
        \caption{The generator $\sigma_{i}$ and $\sigma_{i}^{-1}$.}
        \label{sigma}
    \end{figure}

In this paper, we give a partial refinement of these results. It is well known that any oriented link in the $3$-sphere $S^{3}$ is represented as the closure of a braid, and any braid is presented by a braid word with the generators $\sigma_{i}^{\pm 1}$, which are shown in Figure~\ref{sigma}. See \cite[Chapter 1]{Lickorish-KT}, for example. We denote the set of type $\mathrm{I\hspace{-1.2pt}I^{2}}$ singular fibers of $f$ by $\mathrm{I\hspace{-1.2pt}I^{2}}(f)$, and we denote the cardinality of a finite set $S$ by $| S |$.

    \begin{theorem}
        \label{main1}
            For any link $L$ in $S^{3}$, there exists a stable map $f : S^{3} \to \mathbb{R}^{2}$ without cusp points such that $f$ has no singular fibers of type $\mathrm{I\hspace{-1.2pt}I^{3}}$ and $S_{0}(f)$ is isotopic to $L$. Moreover, if $L$ is represented as the closure of an $n$-braid $\sigma_{z_{1}}^{m_{1}} \sigma_{z_{2}}^{m_{2}} \cdots \sigma_{z_{l}}^{m_{l}}$ with $n \ge 2$,  $1 \leq z_{i} \leq n-1$, $z_{i-1} \neq z_{i}$ and $m_{i}$ non-zero integers, then $|\mathrm{I\hspace{-1.2pt}I^{2}}(f)| = 2(l- | X | )$ holds with $X = \{ i \mid z_{i}=1,~ 1\le i \le l \}$.
            \label{1}
    \end{theorem}

\begin{remark}
While Ishikawa and Koda \cite{Ishikawa-Koda} proved their theorems using $4$-dimensional topology, we establish the theorem above using only $3$-dimensional techniques. Our construction is similar to the method developed in \cite{Ichihara-K}.
\end{remark}

It is well known that every closed orientable $3$-manifold can be obtained from $S^{3}$ by integral Dehn surgery on a link in $S^{3}$. See Section~\ref{s3} for further details. Together with this, we have the following as a corollary of Theorem~\ref{main1}.

\begin{theorem}
    \label{main2}
    Let $M$ be a closed orientable $3$-manifold. Then, there exists a stable map $f_{0}$ from $M$ into $S^{2}$ without cusp points, definite fold points and singular fibers of type $\mathrm{I\hspace{-1.2pt}I^{3}}$. In particular, if $M$ is obtained by an integral Dehn surgery on a link $L$ in $S^{3}$ which is represented as the closure of a pure braid $b=\sigma_{z_{1}}^{m_{1}} \sigma_{z_{2}}^{m_{2}} \cdots \sigma_{z_{i}}^{m_{i}} \cdots \sigma_{z_{l}}^{m_{l}}$, then $|\mathrm{I\hspace{-1.2pt}I^{2}}(f)| = 2(l- | X | )$ holds with $X = \{ i \mid z_{i}=1,~ 1\le i \le l \}$.
\end{theorem}

We remark that, under different assumptions and conditions, there are other studies concerning the construction of stable maps on $3$-manifolds obtained by Dehn surgery. See \cite[Theorem 1.2]{Kalmar-Stipsicz}, for example.

% For example, in \cite[Corollary 4.3]{Ishikawa-Koda}, they showed that Let $M$ be a closed orientable smooth $3$-manifold obtained from $S^{3}$ by Dehn surgery on a non-trivial link $L$ in $S^{3}$. Then there exists a stbale map $f : M \to \mathbb{R}^{2}$ without cusps such that $|\mathrm{I\hspace{-1.2pt}I^{2}} (f)|) \le cr(L)-2$ and $|\mathrm{I\hspace{-1.2pt}I^{3}} (f)|)=\emptyset$, where $cr(L)$ is the crossing number of $L$. The results are refinements of Saeki's result regarding the singularities.

%%%%%%%%%%%%%%%%%%%%%%%%%%%%%%%%%%%%%%%%%%%%%%%%%%%%%%%%%%%%%%%%%%%%%%

\section{A construction of stable maps from $S^{3}$ into the plane}
\label{s2}

In this section, we give a proof of Theorem~\ref{main1}. Before proceeding to the proof, we introduce several definitions and known facts.

For smooth manifolds $M$ and $N$, let $C^{\infty}(M,N)$ be the set of smooth maps from $M$ to $N$ with the Whitney topology. A smooth map $f : M \to N$ is called a \textit{stable map} if there exists a neighborhood $U_{f}$ of $f$ in $C^{\infty}(M,N)$ such that, for any map $g$ in $U_{f}$, there are diffeomorphisms $\Phi : M \to M$ and $\phi: N \to N$ satisfying $g = \phi \circ f \circ \Phi^{-1}$.

In the case where the source and target dimensions are $3$ and $2$, a characterization of stable maps is given as follows (\cite{Levin'65}). For a $3$-manifold $M$, a smooth map $f : M \to \mathbb{R}^{2}$ is a stable map if and only if $f$ is locally described as in one of the following forms:
\begin{enumerate}
    \item $(u,x,y) \mapsto (u,x)$,
    \item $(u,x,y) \mapsto (u,x^{2} + y^{2})$,
    \item $(u,x,y) \mapsto (u,x^{2} - y^{2})$,
    \item $(u,x,y) \mapsto (u,y^{2} + ux - x^{3})$,
\end{enumerate}
%
% \noindent and $f$ globally satisfies
and $f$ globally satisfies

\begin{enumerate}[resume]
    \item $f^{-1} (f(p)) \cap S(f) = \{p\}$ for a cusp point $p$, 
    \item except for cusp points, the restriction of $f$ to $S(f)$ is an immersion with only normal crossings.
\end{enumerate}

The points around which $f$ is described as (2), (3), and (4) are the singular points of $f$, and they are called \textit{definite fold} points, \textit{indefinite fold} points, and \textit{cusp} points, respectively. The set of definite fold points and indefinite fold points is denoted by $S_{0}(f)$ and $S_{1}(f)$, respectively. Consequently, to prove Theorem~\ref{main1}, it suffices to construct a smooth map satisfying these conditions with the described properties. Throughout the following, we denote the boundary of a manifold $N$ by $\partial N$.

\begin{proof}[Proof of Theorem~\ref{main1}.]

Let $L$ be a link in $S^{3}$. Suppose that $L$ is represented by the closure of a braid $b$ of $n$ strings $(n \ge 2)$ given as a braid word $\sigma_{z_{1}}^{m_{1}} \sigma_{z_{2}}^{m_{2}} \cdots \sigma_{z_{l}}^{m_{l}}$, where $1 \leq z_{i} \leq n-1$, $z_{i-1} \neq z_{i}$, $m_{i}$ is non-zero integers for $i \in \{ 1, 2, \ldots, l \}$, and $l \ge 1$.

We will obtain a stable map $S^{3} \to \mathbb{R}^{2}$ by gluing several smooth maps on submanifolds of $S^{3}$. To construct these smooth maps, we begin with the following setup. We decompose $S^{3}$ into two solid tori $V_{1}$ and $V_{2}$ as $S^{3} = V_{1} \cup_{\phi} V_{2}$. Here, $\phi$ denotes a diffeomorphism from $\partial V_{1}$ onto $\partial V_{2}$ such that a meridian of $V_{1}$ is mapped to a longitude of $V_{2}$ by $\phi$. (A $\textit{meridian}$ of $V_{1}$ is $\partial \mathbb{D}^{2} \times \{*\} \subset \mathbb{D}^{2} \times S^{1}=V_{1}$ and a $\textit{longitude}$ of $V_{2}$ is $\{ * \} \times S^{1} \subset \partial \mathbb{D}^{2} \times S^{1} \subset \mathbb{D}^{2} \times S^{1} = V_{2}$.) We further decompose $V_{1}$ into $N_{1}$ and $N_{2}$, each of which is diffeomorphic to $\mathbb{D}^{2} \times [0,l]$. We isotope $L$ so that $L$ is contained in $V_{1}$ and $L \cap N_{1}$ corresponds to the braid $b$. Next, we decompose $N_{1}$ into $l$ solid cylinders $W_{1}, \ldots, W_{l}$, where $W_{i}$ is regarded as $\mathbb{D}^{2} \times [i-1,i]$ $(1 \leq i \leq l)$ such that $L \cap W_{i}$ corresponds to $\sigma_{z_{i}}^{m_{i}}$ as shown in Figure~\ref{V1}. Let $A$ be an annular region in $\mathbb{R}^{2}$ and $E$ a disk in $\mathbb{R}^{2}$ such that the intersection $A \cap E = \partial E$. Also, we decompose $A$ into two rectangular regions $R_{1}$ and $R_{2}$ by $\gamma_{1}, \gamma_{l+1}$. Moreover, we decompose $R_{1}$ into $l$ rectangular regions by arcs $\gamma_{2}, \ldots, \gamma_{l}$ as shown in Figure~\ref{A}.

    \begin{figure}[htbp]
        {\unitlength=1cm
    \begin{picture}(10,7)(0,0)
        \put(0,0){\includegraphics[height=7cm,clip]{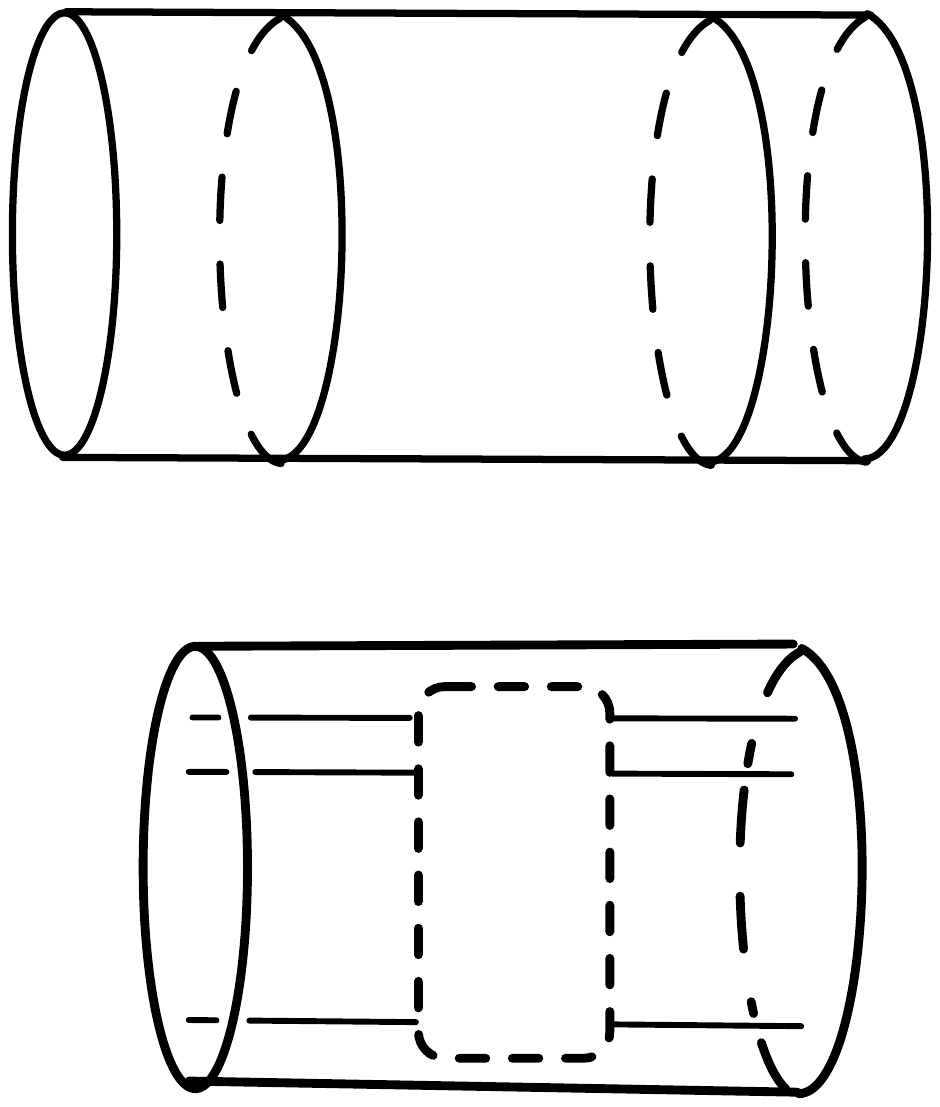}}

        % \put(-0.5,6){$V_{1} :$}
        % \put(-0.5,6.7){$N_{2}$}
        \put(1.5,5.5){$N_{1} :$}

        \put(3.5,3.7){$W_{1}$}
            \put(3,7){$F_{1}$}
            \put(4,7){$F_{2}$}
        \put(5,5.5){$\cdots$}
        \put(7,3.7){$W_{l}$}
            \put(6.5,7){$F_{l}$}
            \put(7.5,7){$F_{l+1}$}

        \put(1.5,2){$W_{i} :$}
        \put(5.3,1.8){$\sigma_{z_{i}}^{m_{i}}$}
        \put(3,3){$1$}
        \put(3,2.5){$2$}
        \put(4.3,1.7){$\vdots$}
        \put(3,1){$n$}
            \put(3.5,0){$F_{i}$}
            \put(7,0){$F_{i+1}$}
        
     \end{picture}}
     \caption{$N_{1} = W_{1} \cup \cdots \cup W_{l}$.}
        \label{V1}
     \end{figure}

\begin{figure}[htbp]
        {\unitlength=1cm
    \begin{picture}(10,6)(0,1)
         \put(2,0){\includegraphics[height=8cm,clip]{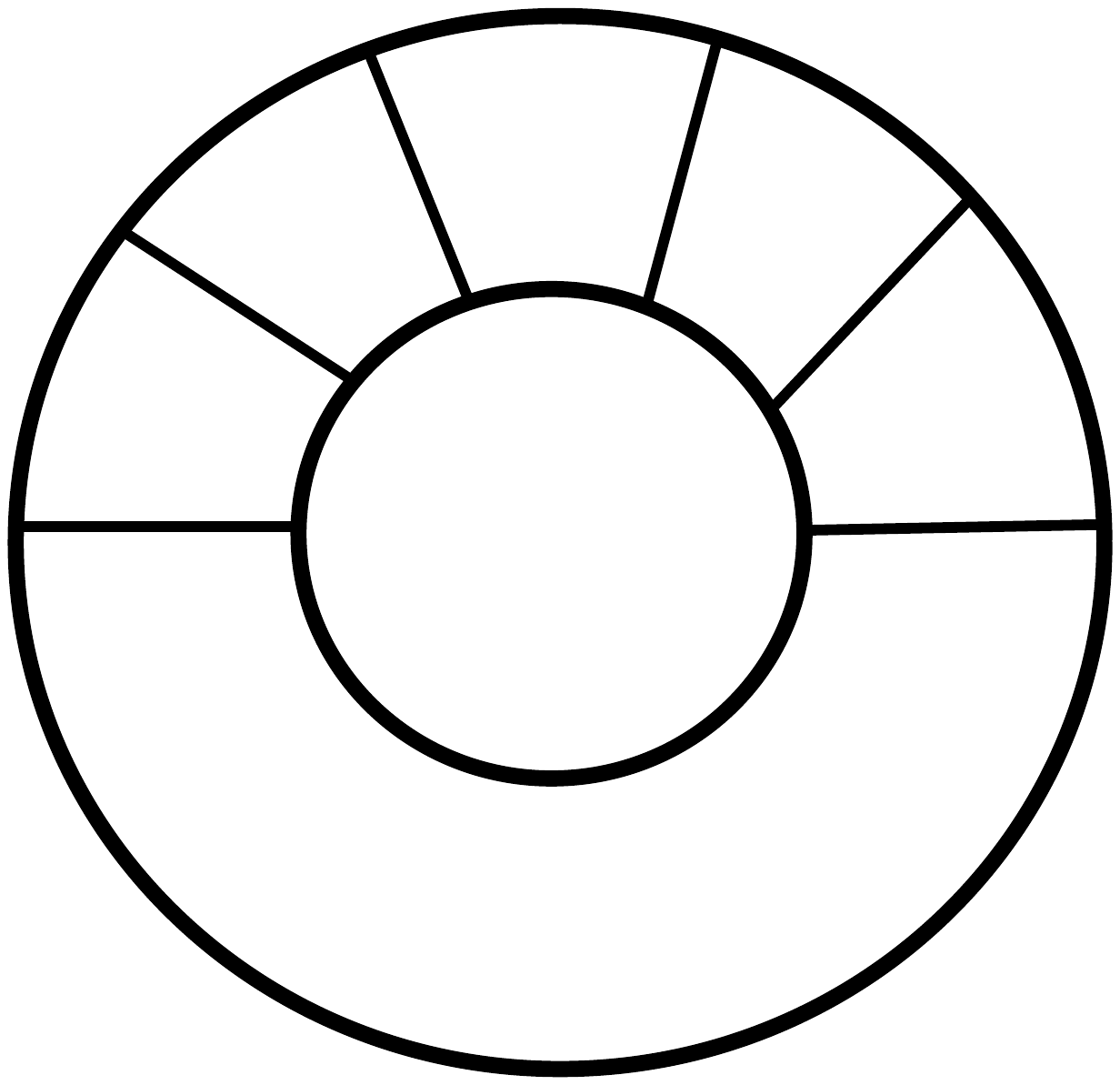}}

    \put(4.7,2.2){$R_{2}$}
         \put(2.7,4.5){$T_{1}$}
            \put(3.5,5.5){$\cdots$}
         \put(4.7,6){$T_{i}$}
            \put(5.7,5.5){$\cdots$}
         \put(6.6,4.5){$T_{l}$}

    \put(1.5,4){$\gamma_{1}$}
    \put(2,5.5){$\gamma_{2}$}
    \put(3.8,6.7){$\gamma_{i}$}
    \put(5.7,6.7){$\gamma_{i+1}$}
    \put(7.3,5.5){$\gamma_{l}$}
    \put(8,4){$\gamma_{l+1}$}

    \end{picture}}
     \caption{$A = R_{1} \cup R_{2}$ and $R_{1} = T_{1} \cup \cdots \cup T_{l}$.}
        \label{A}
\end{figure}

Let $F_{i}$ be the disk $\mathbb{D}^{2} \times \{i-1\}$ in $N_{1}=W_{1} \cup \cdots \cup W_{l}$ $(1 \le i \le l+1)$. Note that $F_{i} = W_{i-1} \cap W_{i}$ $(2 \le i \le l)$. For each $i$, let $\psi_{i} : F_{i} \to \gamma_{i}$ be the map which is obtained by regarding $F_{i}$ as $\mathbb{D}^{2}$ and considering a smooth height function on $\mathbb{D}^{2}$ shown in Figure~\ref{tegaki2}. Precisely, $\psi_{i}$ is the composite map of a diffeomorphism $F_{i} \to \mathbb{D}^{2}$ and a smooth height function on $\mathbb{D}^{2}$ shown in Figure~\ref{tegaki2}. In Figure~\ref{tegaki2}, the fat points on $F_{i}$ depict the local maxima of $\psi_{i}$. For the images of saddle points of $\psi_{i}$, their preimages on $F_{i}$ are figure-eight fibers under $\psi_{i}$. The Reeb graph $\tau_{i}$ is shown in the center of Figure~\ref{tegaki2}. (Let $f$ be a smooth function on a smooth manifold $M$. The \textit{Reeb graph} of $f$ is defined as the quotient space obtained by contracting each connected component of the level sets of $f$ to a point, which often has the structure of a graph.) The boundary $\partial F_i$ is mapped to an endpoint of $\tau_i$.

    \begin{figure}[htbp]
        \setlength\unitlength{1truecm}
        \begin{picture}(15,8)(0,0)
            \put(0,0){\includegraphics[width=1\textwidth,clip]{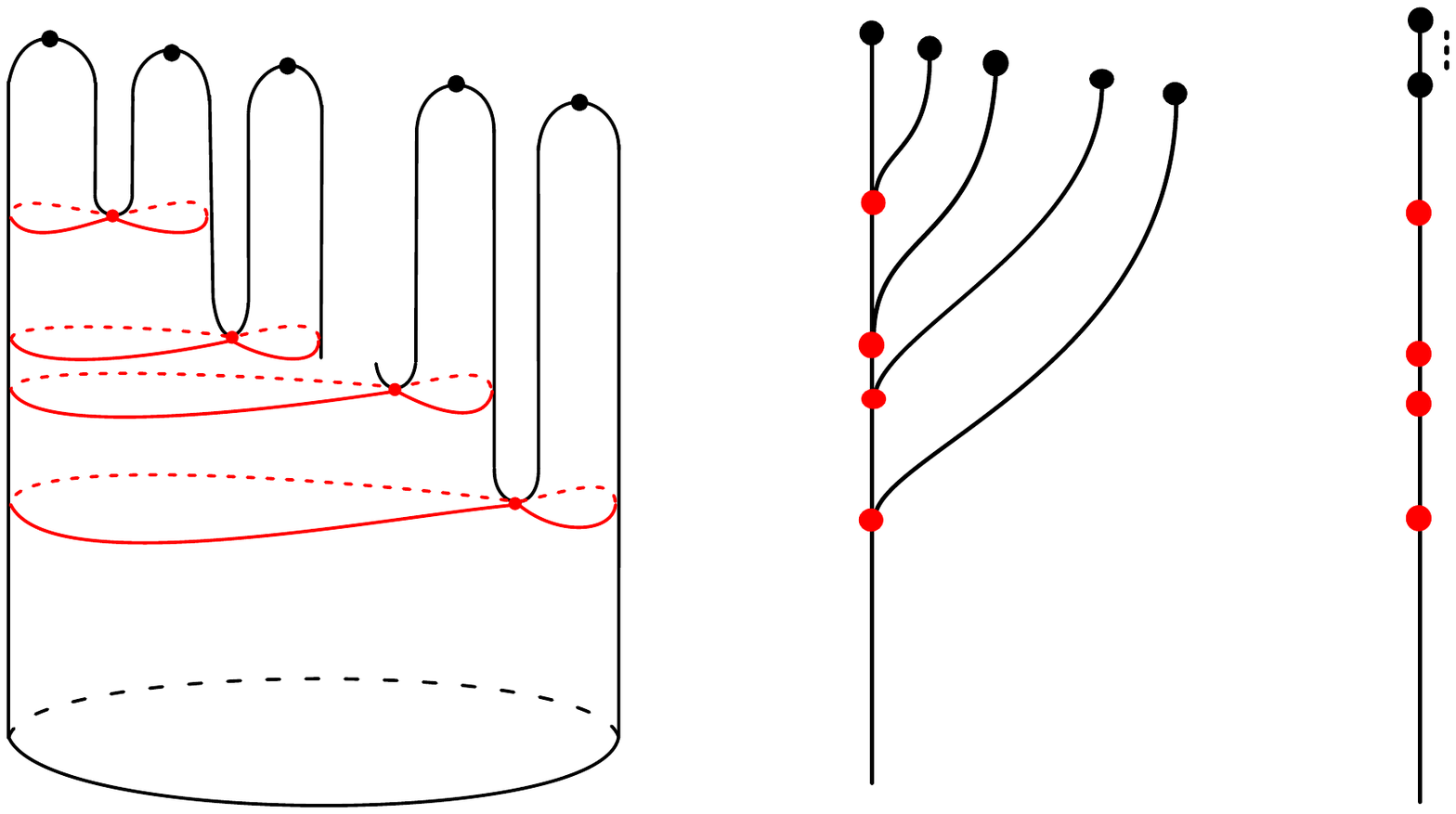}}
            \put(3,6){$\cdots$}
            \put(8.8,7.6){$\cdots$}

        \put(3,0.7){$\mathbb{D}^{2}$}
            \put(0.4,7.9){$1$}
            \put(1.4,7.9){$2$}
            \put(2.4,7.9){$3$}
            \put(3.5,7.9){$n-1$}
            \put(4.9,7.9){$n$}

                \put(6.5,4){$\longrightarrow$}
                
        \put(7.5,0.7){$\tau_{i}$}
            \put(7.4,7.9){$1$}
            \put(7.9,7.9){$2$}
            \put(8.5,7.9){$3$}
            \put(9,7.9){$n-1$}
            \put(10,7.9){$n$}
            
                \put(10.5,4){$\longrightarrow$}

        \put(12,0.7){$\gamma_{i}$}
            \put(12.4,8){$1$}
            \put(12.4,4.6){$\vdots$}
            \put(12.4,7){$n$}

        \put(1,5.7){$P$}
        \put(12.4,6){$\psi_{i}(P)$}
        \put(11.2,6.6){$\alpha$}
        \put(11.5,5.5){ $\stretchleftright{ \{ }{\rule{0ex}{70pt}}{.}$ }
            
        \end{picture}
        \caption{A smooth map $\psi_{i} : F_{i} \to \gamma_{i}$ for $1 \leq i \leq l+1$. The labels $1, 2, \ldots, n$ correspond to the indices of the strands for $\sigma_{z_{i}}^{m_{i}} \subset W_{i}$.}
        \label{tegaki2}
    \end{figure}

For $i$ with $\sigma_{z_{i}}=\sigma_{1}$, we construct a smooth map $\Psi_{i} : W_{i} \to T_{i}$ obtained from an isotopy between $\psi_{i}$ and $\psi_{i+1}$ as follows. Let $P$ be the highest saddle point on $F_{i}$ and $\alpha$ a left-open interval in $\gamma_{i}$ containing $\psi_{i} (P)$ and $\max \mathrm{Im} \psi_{i}$. Then, the isotopy $F_{i} \times [0,1] \to F_{i}$ is given by $m_{i}$ half-twists on the connected component of the preimage of $\alpha$ containing $P$. See Figure~\ref{tegaki1-1} for an illustration of $\Psi_{i}$ when $m_{i}$ is even. Note that the image of definite fold points, i.e., $L \cap W_{i}$, has no normal crossings in this case. When $m_{i}$ is odd, the image of a surface obtained from $F_{i}$ via the isotopy is shown in Figure~\ref{tegaki1-1oF} (left). Further, we deform the height function by another isotopy which is defined by switching the heights of the two leftmost local maxima, as shown in Figure~\ref{tegaki1-1oF}. This induces the switching of the images of the two uppermost maxima of $\mathbb{D}^{2}$ in the target. Thus, when $m_{i}$ is odd, a smooth map $\Psi_{i}$ constructed by connecting the two isotopies is illustrated in Figure~\ref{tegaki1-1o}. Note that the image of definite fold points, i.e., $L \cap W_{i}$, has a normal crossing when $m_{i}$ is odd. Except for this switching, the image of definite fold points has no normal crossings in this case. Furthermore, regardless of the parity of $m_{i}$, the constructed smooth map $\Psi_{i}$ has only definite and indefinite fold points, and has no cusp points.

We make further remarks on Figures~\ref{tegaki1-1} and \ref{tegaki1-1o}. The arcs of $L \cap W_{i}$ are mapped to the segments parallel to the upper edge of $T_{i}$ by $\Psi_{i}$. As noted above, when $m_{i}$ is odd, the segments contain a normal crossing. On $T_{i}$, thick arcs represent parts of the image of indefinite fold points. The boundaries $\partial F_{i}$ and $\partial F_{i+1}$ correspond to the endpoints of the lower edge of $T_{i}$.

\begin{figure}[htbp]
        \setlength\unitlength{1truecm}
        \begin{picture}(15,7)(0,1)
            \put(0,0){\includegraphics[width=1\textwidth,clip]{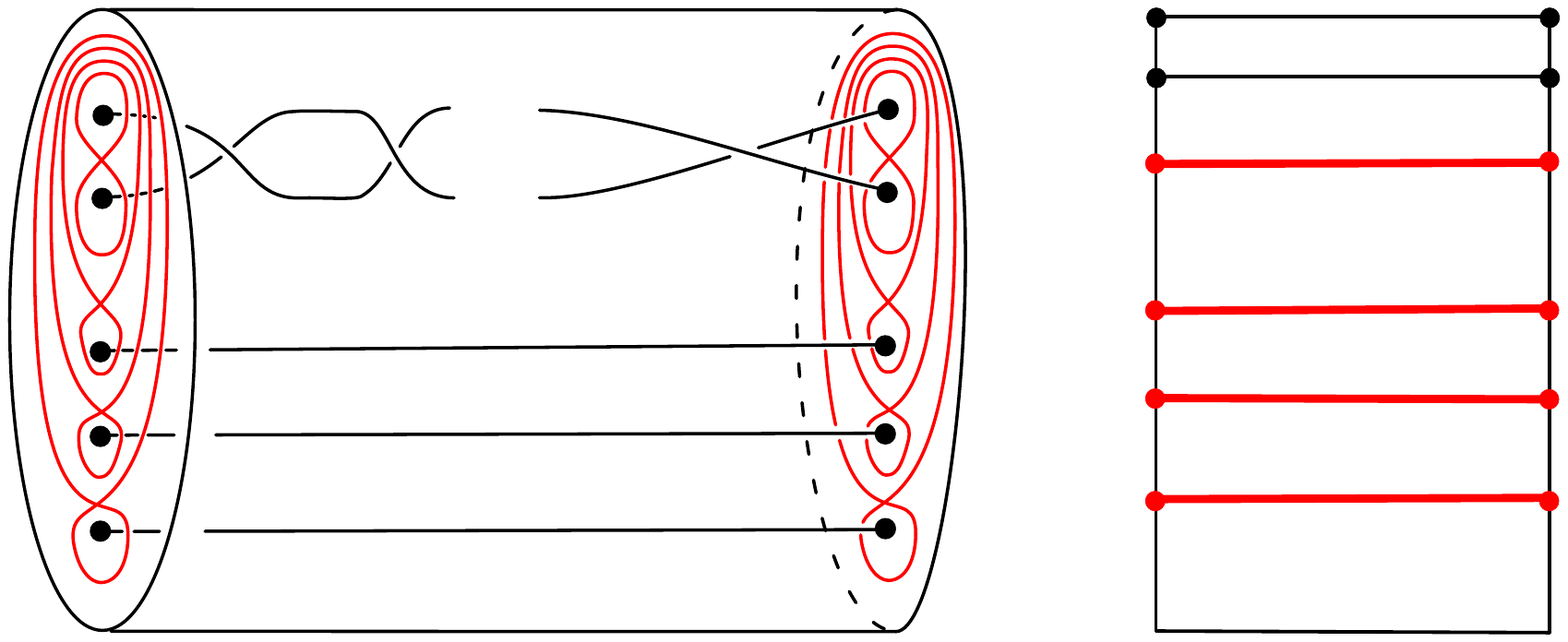}}
                \put(0.5,7.5){$F_{i}$}
                    \put(-0.2,6.2){$1$}
                    \put(-0.2,5.5){$2$}
                    \put(-0.2,4.3){$3$}
                    \put(-0.2,3.6){$4$}
                    \put(-0.2,2.8){$5$}
               
                \put(3.7,5.9){$\cdots$}
                \put(7,7.5){$F_{i+1}$}
                \put(8.8,7.5){$\gamma_{i}$}
                \put(11.9,7.5){$\gamma_{i+1}$}
                    %\put(12.3,7){$1,\cdots,5$}

            \put(8,4.5){$\longrightarrow$}
                \put(10.5,6.65){$\vdots$}
                \put(8.7,6.9){$1$}
                \put(8.7,6.4){$5$}

            \put(4,1.5){$W_{i}$}
            \put(10.5,1.5){$T_{i}$}
        \end{picture}
        \caption{A stable map $\Psi_{i} : W_{i} \to T_{i}$, when $z_{i}=1$, $n=5$, and $m_{i}$ is even.}
        \label{tegaki1-1}
    \end{figure}

\begin{figure}[htbp]
        \setlength\unitlength{1truecm}
        \begin{picture}(15,9)(0,0)
            \put(0,0){\includegraphics[width=1\textwidth,clip]{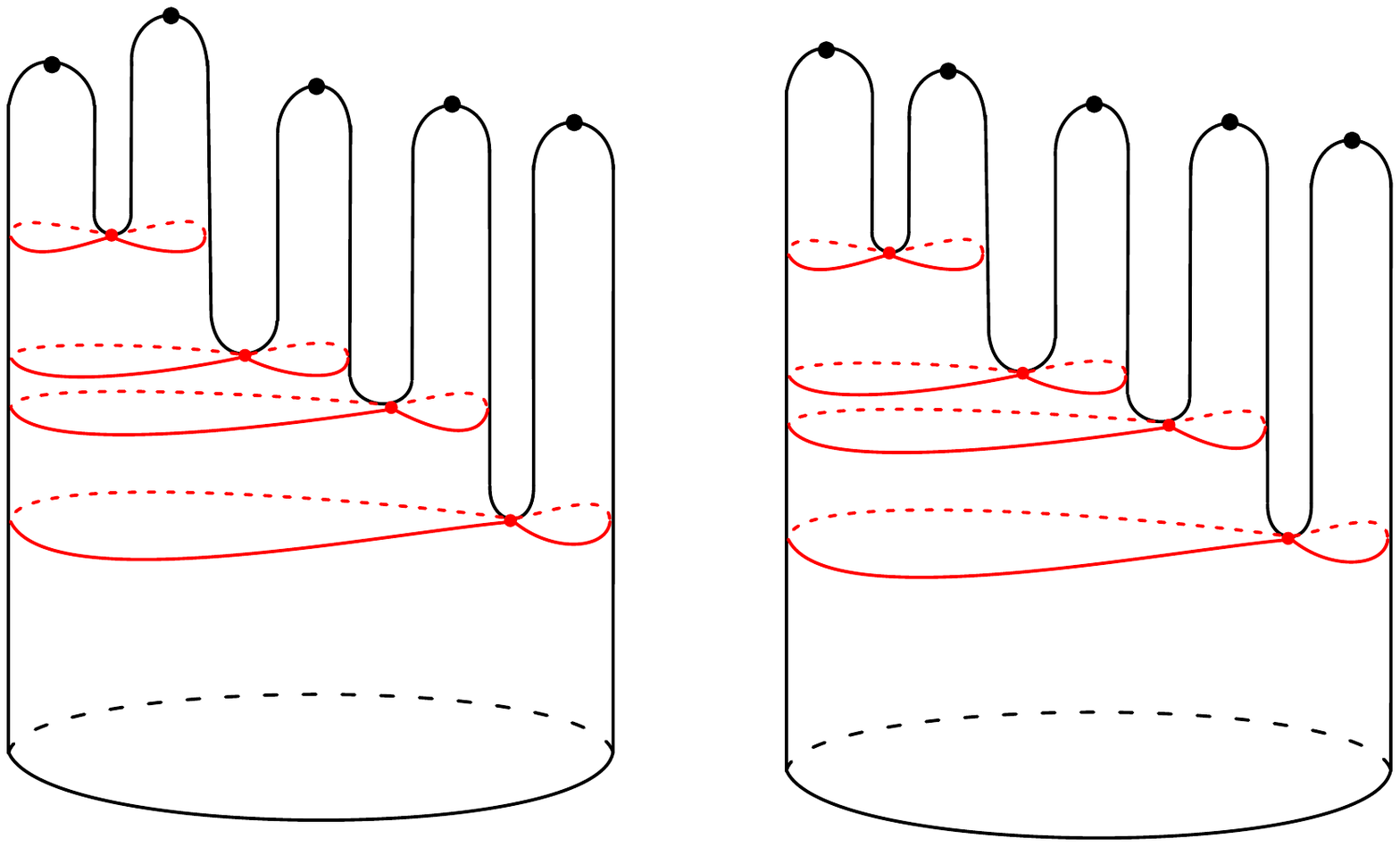}}

        %\put(2.5,0){$(\text{A})$}
            \put(0.4,8){$2$}
            \put(1.4,8.3){$1$}
            \put(2.7,8){$3$}
            \put(3.8,7.9){$4$}
            \put(5,7.7){$5$}

        \put(5.8,5){$\longrightarrow$}
                
        %\put(9,0){$F_{i+1}$}
            \put(7.1,8.3){$2$}
            \put(8.2,8){$1$}
            \put(9.4,7.8){$3$}
            \put(10.5,7.7){$4$}
            \put(11.5,7.5){$5$}
            
        \end{picture}
        \caption{The switching of local maxima.}
        \label{tegaki1-1oF}
    \end{figure}

    \begin{figure}[htbp]
        \setlength\unitlength{1truecm}
        \begin{picture}(15,7)(0,1)
            \put(0,0){\includegraphics[width=1\textwidth,clip]{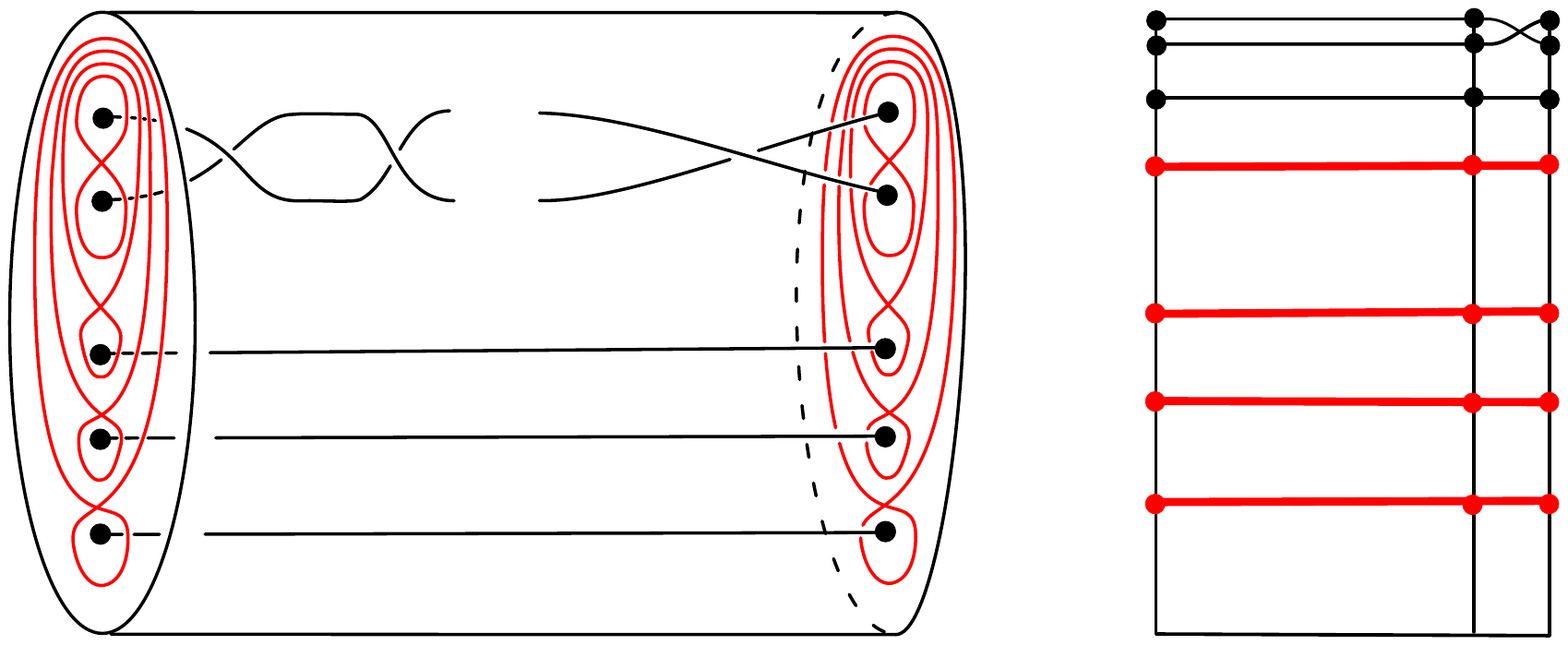}}
                \put(0.5,7.5){$F_{i}$}
                    \put(-0.2,6.2){$1$}
                    \put(-0.2,5.5){$2$}
                    \put(-0.2,4.3){$3$}
                    \put(-0.2,3.6){$4$}
                    \put(-0.2,2.8){$5$}
               
                \put(3.7,5.9){$\cdots$}
                \put(7,7.5){$F_{i+1}$}
                \put(8.8,7.5){$\gamma_{i}$}
                \put(11.9,7.5){$\gamma_{i+1}$}
                % \put(11.2,7.5){$(\text{A})$}
                    %\put(12.3,7){$1,\cdots,5$}

            \put(8,4.5){$\longrightarrow$}
                \put(10.5,6.45){$\vdots$}
                \put(8.7,7){$1$}
                \put(8.7,6.7){$2$}
                \put(8.7,6.3){$5$}
                
            \put(4,1.5){$W_{i}$}
            \put(10.5,1.5){$T_{i}$}
        \end{picture}
        \caption{A stable map $\Psi_{i} : W_{i} \to T_{i}$, when $z_{i}=1$, $n=5$, and $m_{i}$ is odd.}
        \label{tegaki1-1o}
    \end{figure}

For $i$ with $\sigma_{z_{i}} \neq \sigma_{1}$, to construct a smooth map $\Psi_{i} : W_{i} \to T_{i}$, we first consider the deformation between $\psi_{i}$ and $\psi_{i+1}$ illustrated in Figure~\ref{<3}. In Figure~\ref{<3}, the deformation starts at the top corresponding to $\psi_{i}$, goes down to obtain another smooth map $\psi_{i}^{\prime} : F_{i}^{\prime} \to \gamma_{i}^{\prime}$, goes down to obtain another smooth map $\psi_{i}^{\prime \prime} : F_{i}^{\prime \prime} \to \gamma_{i}^{\prime \prime}$, and gets back conversely to the top, corresponding to $\psi_{i+1}$. We denote the intermediate disk (middle) by $F_{i}^{\prime}$ and the intermediate disk (bottom) by $F_{i}^{\prime \prime}$. Note that during the deformation in Figure~\ref{<3}, the singular points and singular fibers move on the disks. The movement can be understood by observing the heights of two adjacent saddle points. At $F_{i}^{\prime}$, the heights of these two adjacent saddle points lie at the same level. At $F_{i}^{\prime \prime}$, the relative positions of these two adjacent saddle points are switched.

We next consider an isotopy from $\psi_{i}^{\prime}$ to $\psi_{i}^{\prime}$ through $\psi_{i}^{\prime \prime}$. The isotopy is determined by $m_{i}$ half-twists in the same way as in the case $z_{i}=1$. When $m_{i}$ is even, we plug in the isotopy during the deformation from $\psi_{i}^{\prime}$ to $\psi_{i}^{\prime}$. Then, a smooth map $\Psi_{i} : W_{i} \to T_{i}$ is obtained by the deformation with isotopy, which is illustrated in Figure~\ref{tegaki1-2}. In Figure~\ref{tegaki1-2}, $F_{i}^{\prime}$ appears twice in $W_{i}$ on both sides of the crossings of $L$ corresponding to $\sigma_{z_{i}}^{m_{i}}$. On each $F_{i}^{\prime}$, there is a singular fiber of type $\mathrm{I\hspace{-1.2pt}I^{2}}$. When $m_{i}$ is odd, we construct the smooth map $\Psi_{i}$ by combining two isotopies. The first is the isotopy during the deformation from $\psi_{i}^{\prime}$ to $\psi_{i}^{\prime}$, and the second is defined by switching the heights of the $z_{i}$-th and $(z_{i}+1)$-th local maxima from the left in Figure~\ref{tegaki1-2oF}. Similarly to the case where $z_{i}=1$, this switching creates a normal crossing in the image of $S_0(f)$. Consequently, as illustrated in Figure~\ref{tegaki1-2oF}, the map $\Psi_{i}$ is almost identical to the case where $m_{i}$ is even, except for the presence of this normal crossing. Regardless of the parity of $m_{i}$, the constructed smooth map $\Psi_{i}$ has definite fold points, indefinite fold points, and only two singular fibers of type $\mathrm{I\hspace{-1.2pt}I^{2}}$, and has no cusp points.

    \begin{figure}[htbp]
        \setlength\unitlength{1truecm}
        \begin{picture}(15,20)(0,0)
            \put(0.5,13){\includegraphics[width=0.8\textwidth,clip]{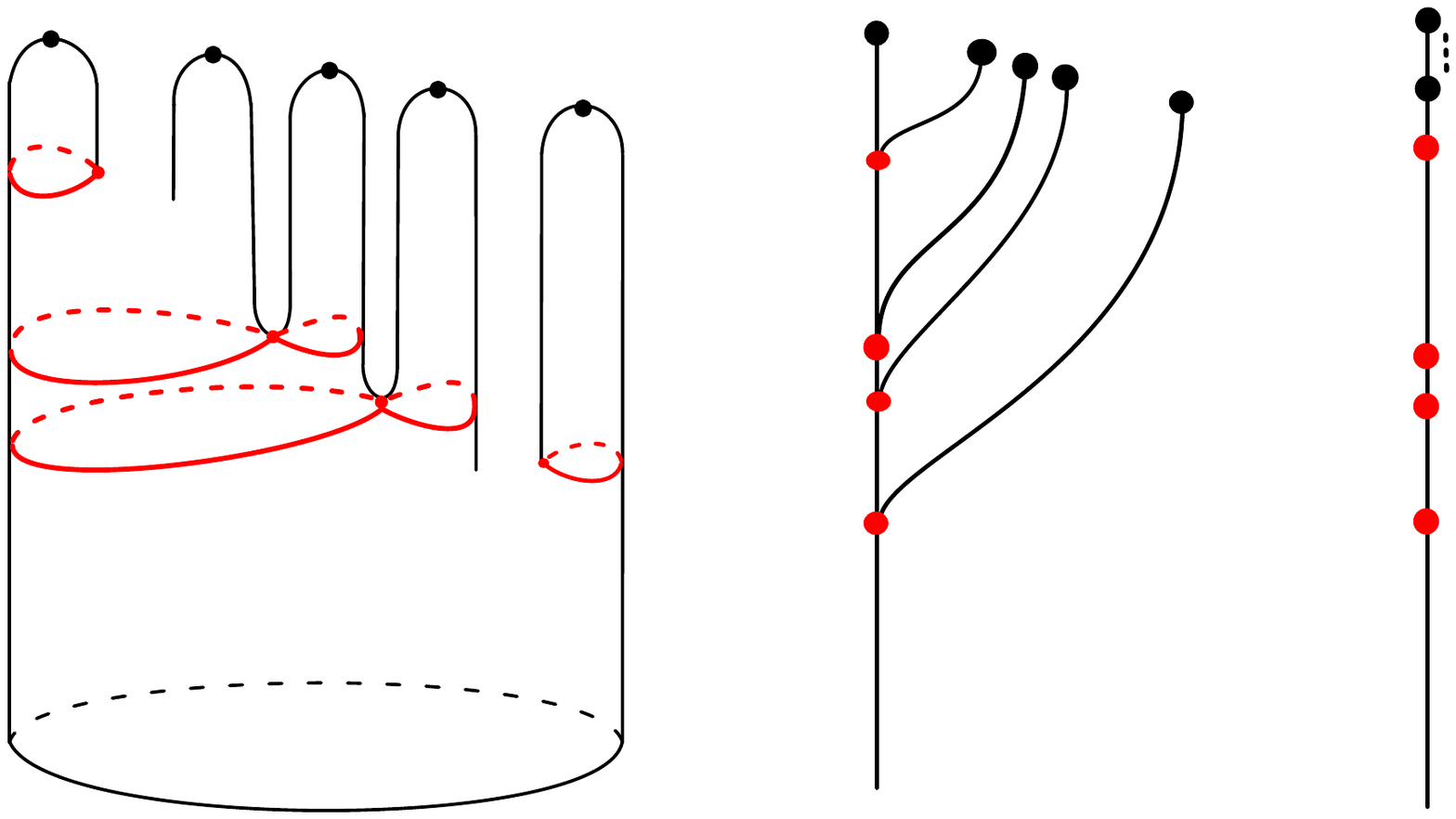}}
            \put(2.3,13.5){$F_{i}=F_{i+1}$}
                \put(0.8,19.3){$1$}
                \put(1.6,19.3){$i-1$}
                \put(2.7,19.3){$i$}
                \put(3.2,19.3){$i+1$}
                \put(4.5,19.3){$l$}
            \put(1.3,18.5){$\cdots$}
            \put(3.8,18){$\cdots$}
                    \put(6.4,19.3){$1$}
                    %\put(6.7,19.5){$i-1$}
                    \put(7.4,19.3){$i$}
                    %\put(7.7,19.5){$i+1$}
                    \put(8.6,19.3){$l$}
            \put(6.7,19){$\cdots$}
            \put(8,18.5){$\cdots$}
                    \put(10.5,19.2){$1$}
                        \put(10.5,17.5){$\vdots$}
                        \put(10.5,16.1){$\vdots$}
                    \put(10.5,18.7){$l$}
                \put(-1.5,16.5){$\psi_{i} = \psi_{i+1} : $}
                \put(5.5,16.5){$\to$}
                \put(9,16.5){$\to$}

            \put(0.5,6.5){\includegraphics[width=0.8\textwidth,clip]{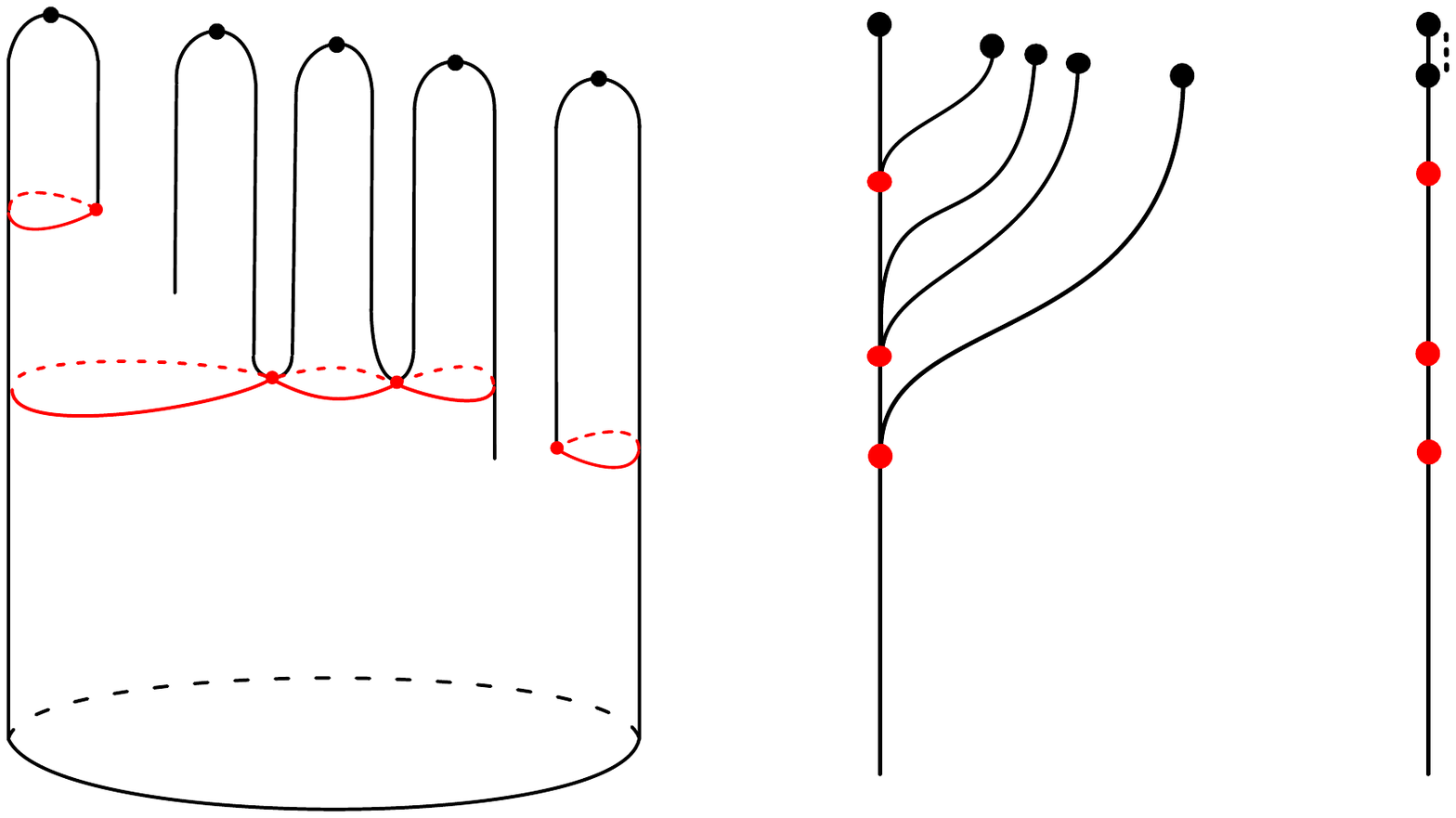}}
            \put(2.7,6.9){$F_{i}^{\prime}$}
                \put(0.8,12.9){$1$}
                \put(1.7,12.9){$i-1$}
                \put(2.8,12.9){$i$}
                \put(3.3,12.9){$i+1$}
                \put(4.6,12.9){$l$}
            \put(1.3,12){$\cdots$}
            \put(3.95,12){$\cdots$}
                    \put(6.4,12.8){$1$}
                    %\put(6.7,13){$i-1$}
                    \put(7.5,12.8){$i$}
                    %\put(7.7,13){$i+1$}
                    \put(8.6,12.8){$l$}
            \put(6.7,12.5){$\cdots$}
            \put(8,12){$\cdots$}
                    \put(10.5,12.6){$1$}
                        \put(10.5,11){$\vdots$}
                        \put(10.5,9.9){$\vdots$}
                    \put(10.5,12){$l$}
                \put(0,10){$\psi_{i}^{\prime} : $}
                \put(5.5,10){$\to$}
                \put(9,10){$\to$}

            \put(0.5,0){\includegraphics[width=0.8\textwidth,clip]{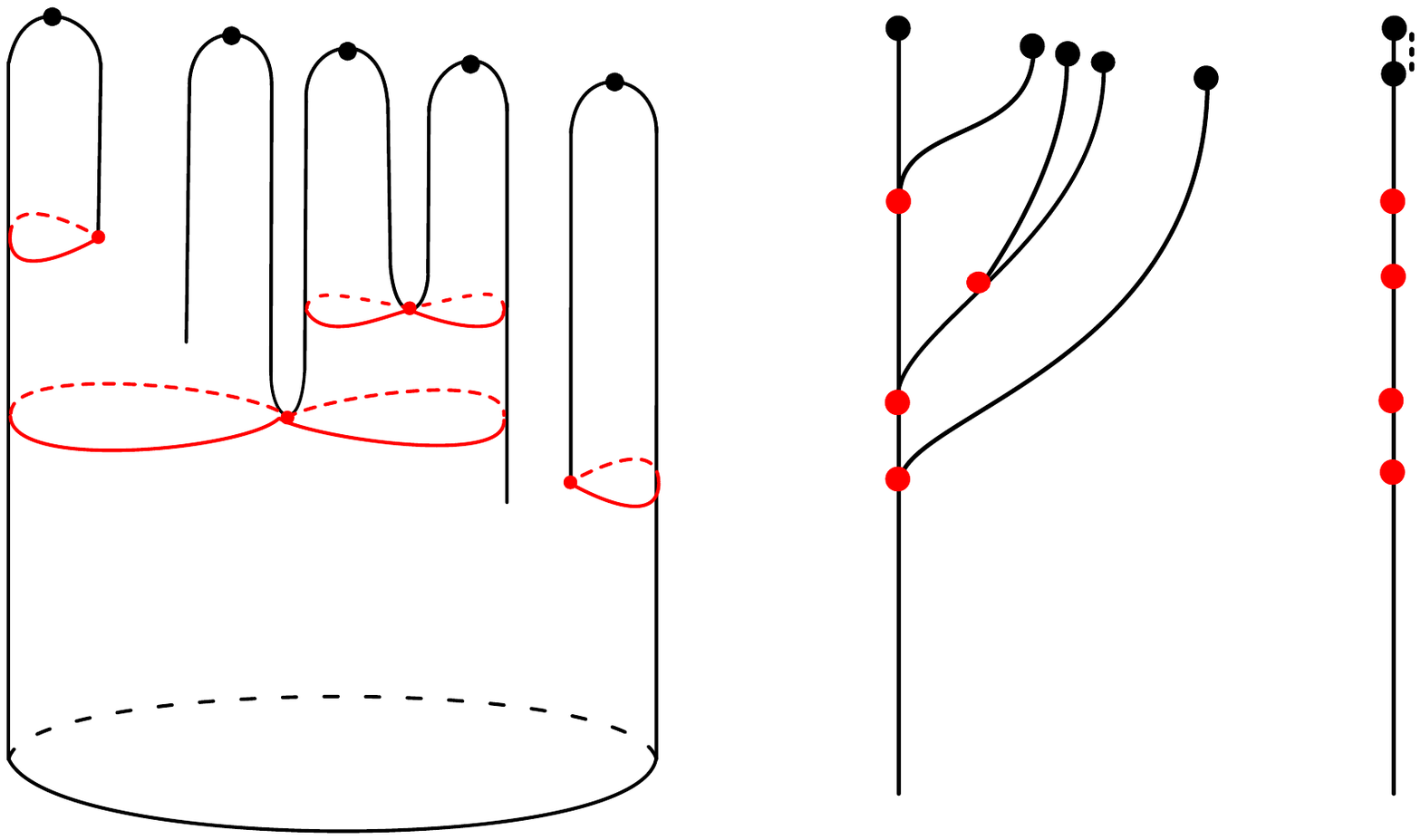}}
            \put(2.7,0.3){$F_{i}^{\prime \prime}$}
                \put(1,6.5){$1$}
                \put(2,6.5){$i-1$}
                \put(3,6.5){$i$}
                \put(3.6,6.5){$i+1$}
                \put(4.9,6.5){$l$}
            \put(1.5,5){$\cdots$}
            \put(4.2,5){$\cdots$}
                    \put(6.8,6.5){$1$}
                    %\put(7.3,6.8){$i-1$}
                    \put(8,6.5){$i$}
                    %\put(8.1,6.8){$i+1$}
                    \put(9.1,6.5){$l$}
            \put(7.2,6.2){$\cdots$}
            \put(8.4,5.8){$\cdots$}
                    \put(10.5,6.3){$1$}
                        \put(10.5,4.7){$\vdots$}
                        \put(10.5,3.3){$\vdots$}
                    \put(10.5,5.8){$l$}
                \put(0,3.6){$\psi_{i}^{\prime \prime} : $}
                \put(5.5,3.6){$\to$}
                \put(9,3.6){$\to$}
            
        \end{picture}
        \caption{A deformation from $\psi_{i}$ to $\psi_{i+1}$, when $z_{i} \ge 2$.}
        \label{<3}
    \end{figure}

    \begin{figure}[htbp]
        \setlength\unitlength{1truecm}
        \begin{picture}(15,6.5)(0,1)
            \put(0,0){\includegraphics[width=1\textwidth,clip]{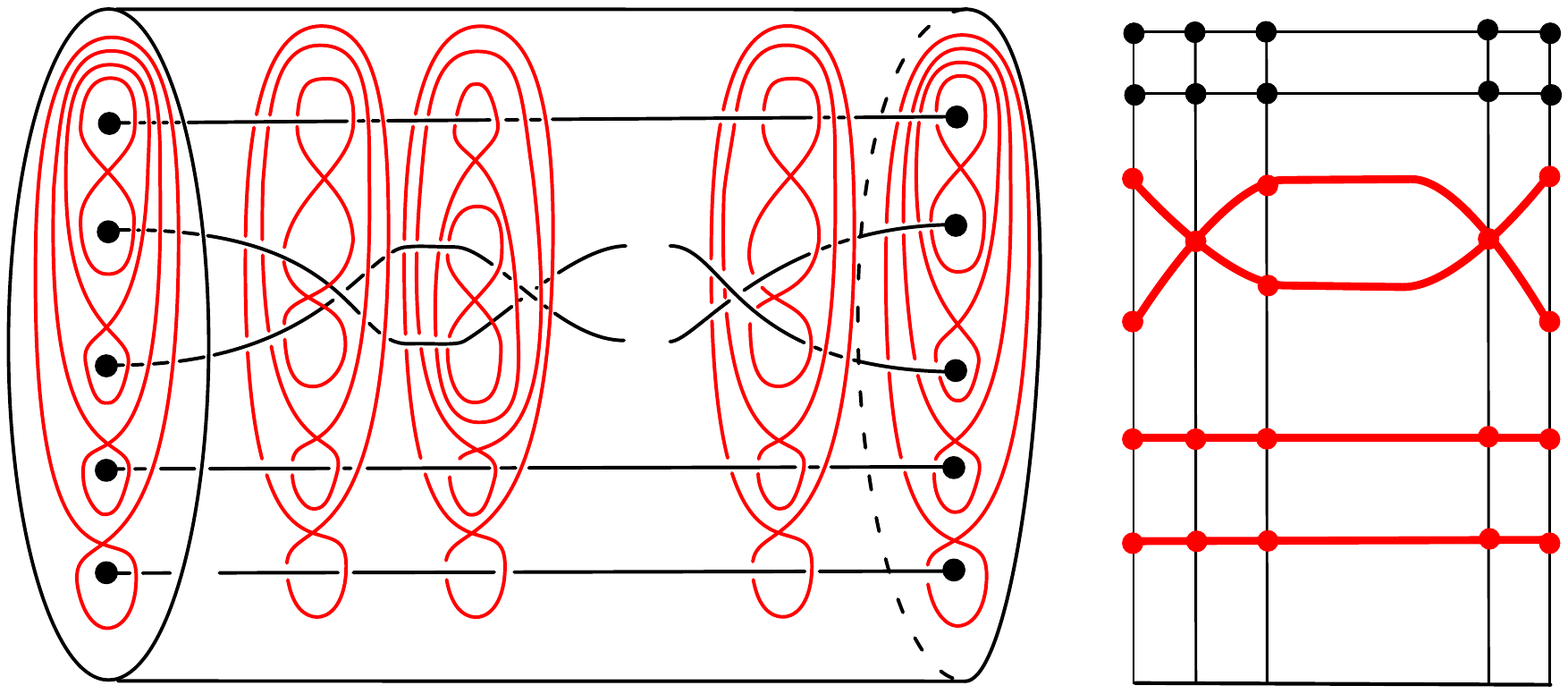}}
            \put(0.5,7.5){$F_{i}$}
                %\put(2.5,7.5){$\overbrace{}^{m_{i}}$}
                \put(7.5,7.5){$F_{i+1}$}
                    % \put(-0.2,6.2){$1$}
                    % \put(-0.2,5.5){$2$}
                    % \put(-0.2,4.3){$3$}
                    % \put(-0.2,3.6){$4$}
                    % \put(-0.2,2.8){$5$}
                
                \put(8.9,7.5){$\gamma_{i}$}
                \put(9.5,7.5){$\gamma_{i}^{\prime}$}
                \put(10,7.5){$\gamma_{i}^{\prime \prime}$}
                \put(11.8,7.5){$\gamma_{i}^{\prime}$}
                \put(12.2,7.5){$\gamma_{i+1}$}
                    \put(11,6.75){$\vdots$}

                \put(2.4,7.5){$F_{i}^{\prime}$}
                \put(3.5,7.5){$F_{i}^{\prime \prime}$}
                \put(5,5){$\cdots$}
                \put(6.2,7.5){$F_{i}^{\prime}$}
                    % \put(13,7){$1,\cdots,5$}

            \put(8.5,4.5){$\longrightarrow$}

            \put(4,1.5){$W_{i}$}
            \put(10.5,1.5){$T_{i}$}
        \end{picture}
        \caption{When $n=5$, $z_{i}=2$ and $m_{i}$ is even, a stable map $\Psi_{i} : W_{i} \to T_{i}$.}
        \label{tegaki1-2}
    \end{figure}

    \begin{figure}[htbp]
        \setlength\unitlength{1truecm}
        \begin{picture}(15,9)(0,0)
            \put(0,0){\includegraphics[width=1\textwidth,clip]{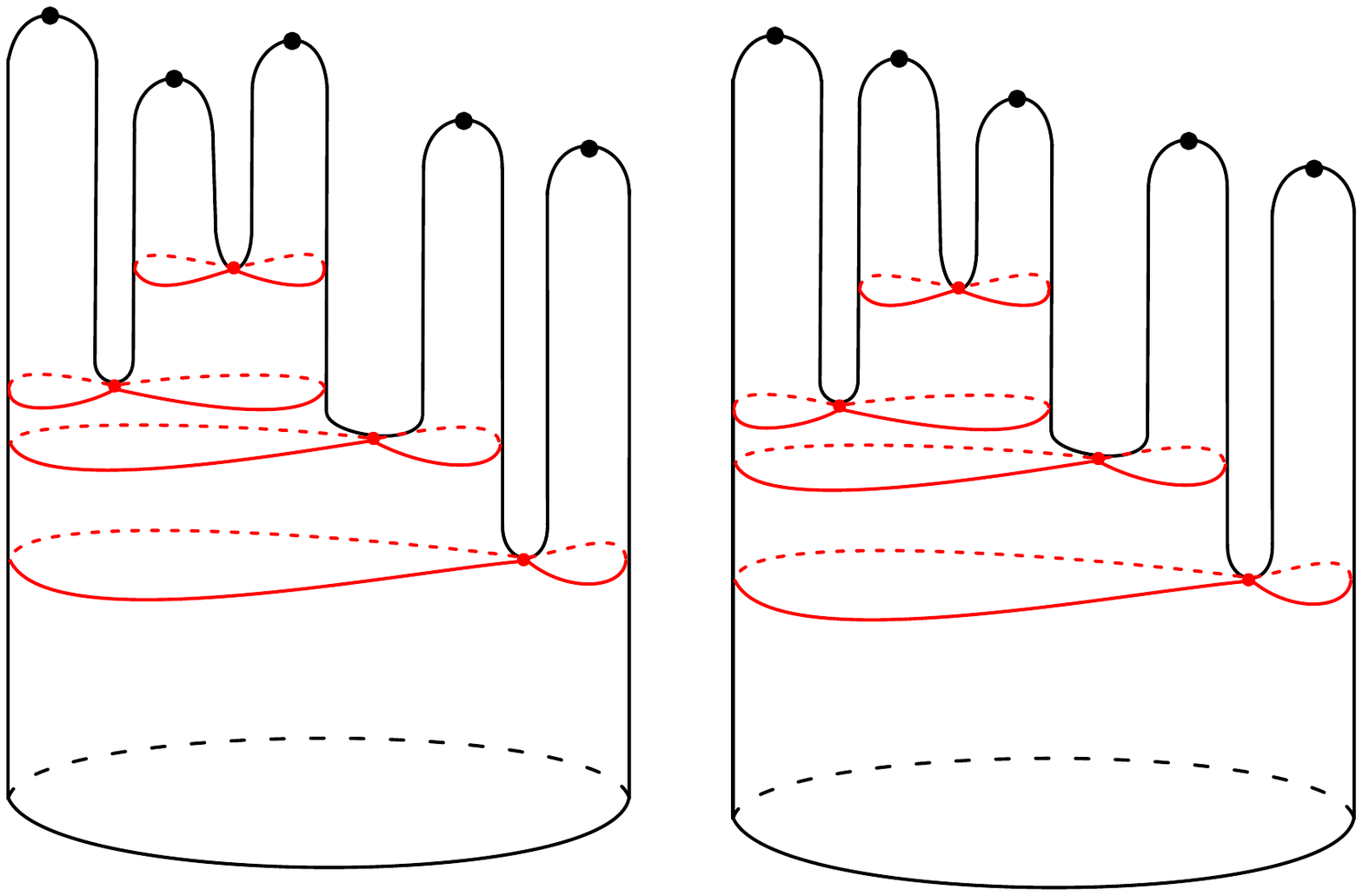}}

        %\put(2.5,0){$(\text{A})$}
            \put(0.6,8.8){$1$}
            \put(1.7,8.3){$3$}
            \put(2.8,8.6){$2$}
            \put(4.4,7.9){$4$}
            \put(5.5,7.7){$5$}

        \put(6,5){$\longrightarrow$}
                
        %\put(9,0){$F_{i+1}$}
            \put(7.1,8.8){$1$}
            \put(8.2,8.5){$3$}
            \put(9.4,8){$2$}
            \put(10.8,7.7){$4$}
            \put(11.8,7.5){$5$}
            
        \end{picture}
        \caption{The switching of local maxima when $z_{i}=2$.}
        \label{tegaki1-2oF}
    \end{figure}

    \begin{figure}[htbp]
        \setlength\unitlength{1truecm}
        \begin{picture}(15,6.5)(0,1)
            \put(0,0){\includegraphics[width=1\textwidth,clip]{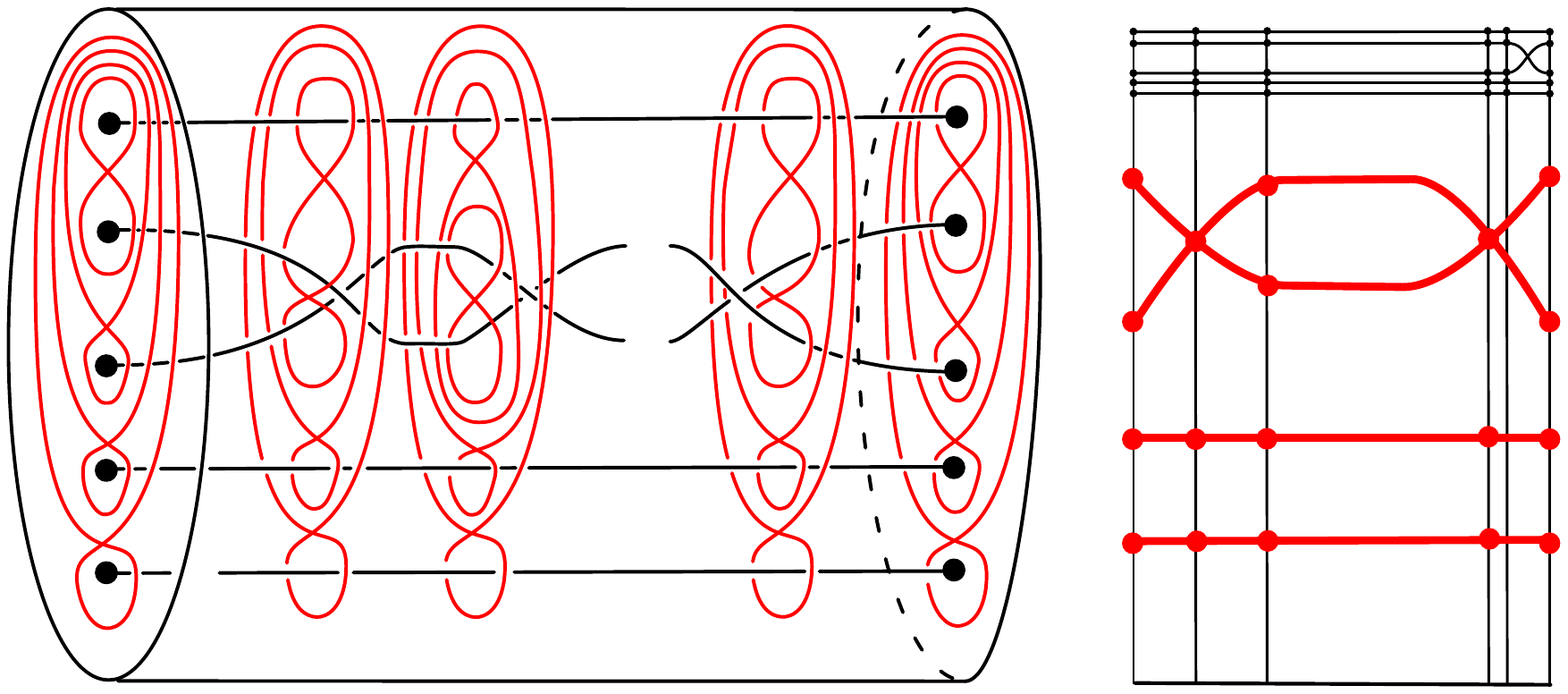}}
            \put(0.5,7.5){$F_{i}$}
                %\put(2.5,7.5){$\overbrace{}^{m_{i}}$}
                \put(7.5,7.5){$F_{i+1}$}
                    % \put(-0.2,6.2){$1$}
                    % \put(-0.2,5.5){$2$}
                    % \put(-0.2,4.3){$3$}
                    % \put(-0.2,3.6){$4$}
                    % \put(-0.2,2.8){$5$}
                
                \put(8.9,7.5){$\gamma_{i}$}
                \put(9.5,7.5){$\gamma_{i}^{\prime}$}
                \put(10,7.5){$\gamma_{i}^{\prime \prime}$}
                \put(11.8,7.5){$\gamma_{i}^{\prime}$}
                \put(12.2,7.5){$\gamma_{i+1}$}
                    \put(11,6.75){$\vdots$}

                \put(2.4,7.5){$F_{i}^{\prime}$}
                \put(3.5,7.5){$F_{i}^{\prime \prime}$}
                \put(5,5){$\cdots$}
                \put(6.2,7.5){$F_{i}^{\prime}$}
                    % \put(13,7){$1,\cdots,5$}

            \put(8.5,4.5){$\longrightarrow$}

            \put(4,1.5){$W_{i}$}
            \put(10.5,1.5){$T_{i}$}
        \end{picture}
        \caption{When $n=5$, $z_{i}=2$ and $m_{i}$ is odd, a stable map $\Psi_{i} : W_{i} \to T_{i}$.}
        \label{tegaki1-2}
    \end{figure}

By connecting $\Psi_{1} : W_{1} \to T_{1}, \ldots, \Psi_{l} : W_{l} \to T_{l}$ along $\psi_{2} : F_{2} \to \gamma_{2}, \ldots, \psi_{l} : F_{l} \to \gamma_{l}$, we obtain a smooth map $N_{1} \to R_{1}$. We connect this map with a smooth map $N_{2} \to R_{2}$ naturally induced from $\psi_{1} : F_{1} \to \gamma_{1}$ and $\psi_{l+1} : F_{l+1} \to \gamma_{l+1}$ by the product structure of $N_{2}$. As a result, a smooth map $g_{1} : V_{1} \to A$ is obtained. Further, connecting $g_{1}$ and the natural projection $g_{2} : V_{2} \to E$, we finally obtain a smooth map $f : S^{3} \to A \cup E \subset \mathbb{R}^{2}$. By construction, this map $f$ has no cusp points, and its singular set consists of definite and indefinite fold points. Also, $f$ satisfies the global conditions (5) and (6). Moreover, the set of definite fold points forms a link isotopic to $L$ and $f$ has no singular fibers of type $\mathrm{I\hspace{-1.2pt}I^{3}}$. Consequently, $f : S^{3} \to \mathbb{R}^{2}$ is a stable map without cusp points such that $f$ has no singular fibers of type $\mathrm{I\hspace{-1.2pt}I^{3}}$ and $S_{0}(f)$ is isotopic to $L$.

Furthermore, $|\mathrm{I\hspace{-1.2pt}I^{2}}(f)| = 2(l-X)$ holds by the construction of $f$, since $L$ is represented as the closure of an $n$-braid $\sigma_{z_{1}}^{m_{1}} \sigma_{z_{2}}^{m_{2}} \cdots \sigma_{z_{l}}^{m_{l}}$ with $n \ge 2$, $1 \leq z_{i} \leq n-1$, $z_{i-1} \neq z_{i}$, $m_{i} \ne 0$.

\end{proof}
%%%%%%%%%%%%%%%%%%%%%%%%%%%%%%%%%%%%%%%%%%%%%

\begin{example}

Following the construction in the proof above, we provide a stable map $f : S^{3} \to \mathbb{R}^{2}$, shown in Figure~\ref{tegaki8}, 
such that $S_{0}(f)$ is isotopic to $K=10_{128}$ in the standard knot table. It is known that $K$ is represented by a $4$-string braid with the braid word $\sigma_1^{-3} \sigma_2^{-1} \sigma_1^{-2} \sigma_2^{-2} \sigma_3^{-1} \sigma_2 \sigma_3^{-1}$. For details, see \cite{KnotInfo}. In Figure~\ref{tegaki8} (right), the image of $S^{3}$ by $f$ is a disk in $\mathbb{R}^{2}$ such that $f ( S_{1}(f))$ is the thick immersed curve in the disk. Also, $S_0(f)$, which is isotopic to $K$, is the immersed curve depicted in Figure~\ref{tegaki8} (left). In Figure~\ref{tegaki8} (right), the box $(\text{B})$ represents a specific crossing pattern as illustrated in Figure~\ref{Box}. This $f$ has no singular fibers of type $\mathrm{I\hspace{-1.2pt}I^{3}}$. Only some of the singular fibers are shown in the figure.

    \begin{figure}[htbp]
        \setlength\unitlength{1truecm}
        \begin{picture}(15,11)(0,-1)
            \put(-0.3,0){\includegraphics[width=0.5\textwidth,clip]{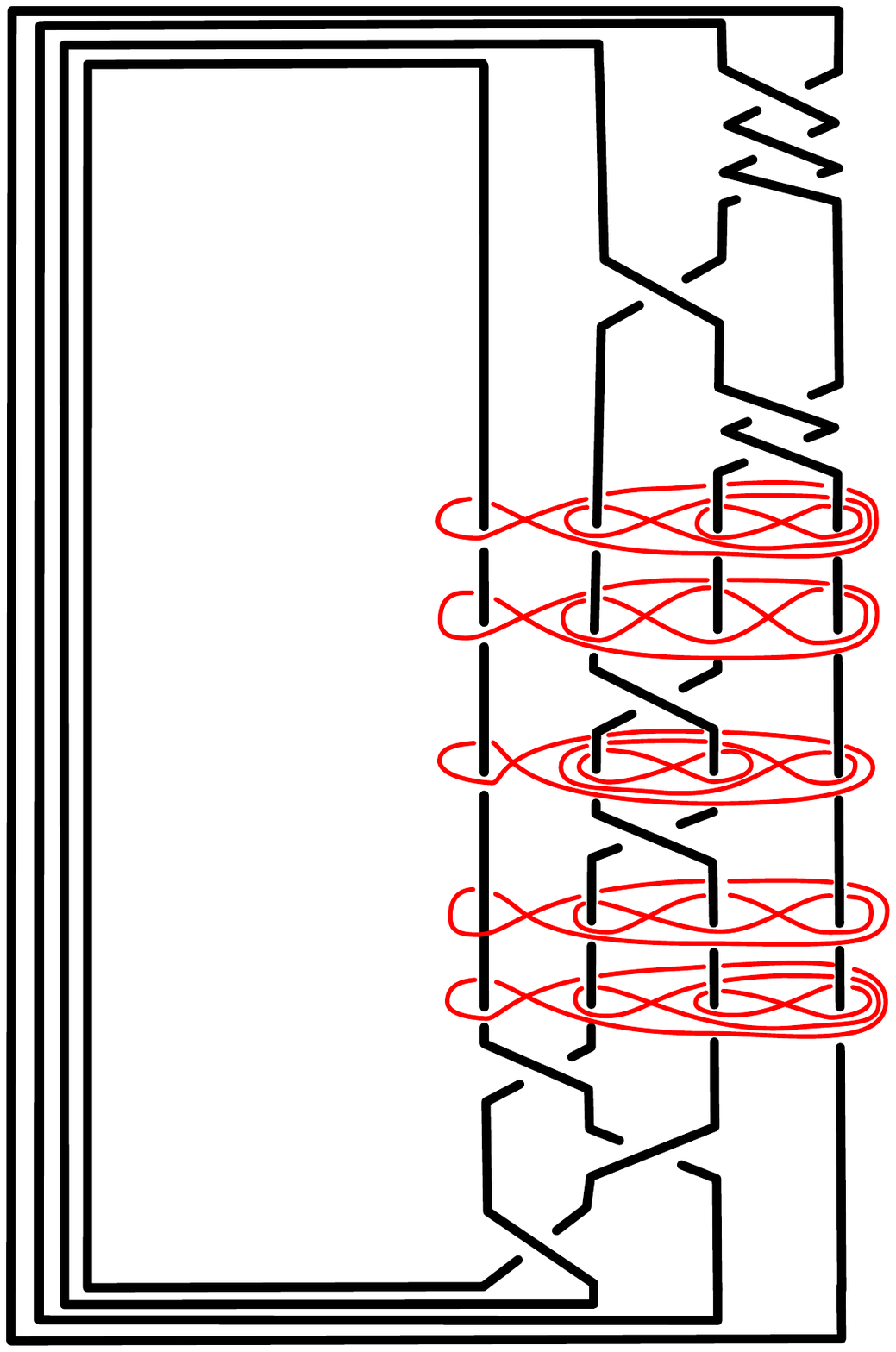}}
            \put(3,-0.3){$K$}

            \put(6,4.5){$\longrightarrow$}

            \put(6.6,0){\includegraphics[width=0.5\textwidth,clip]{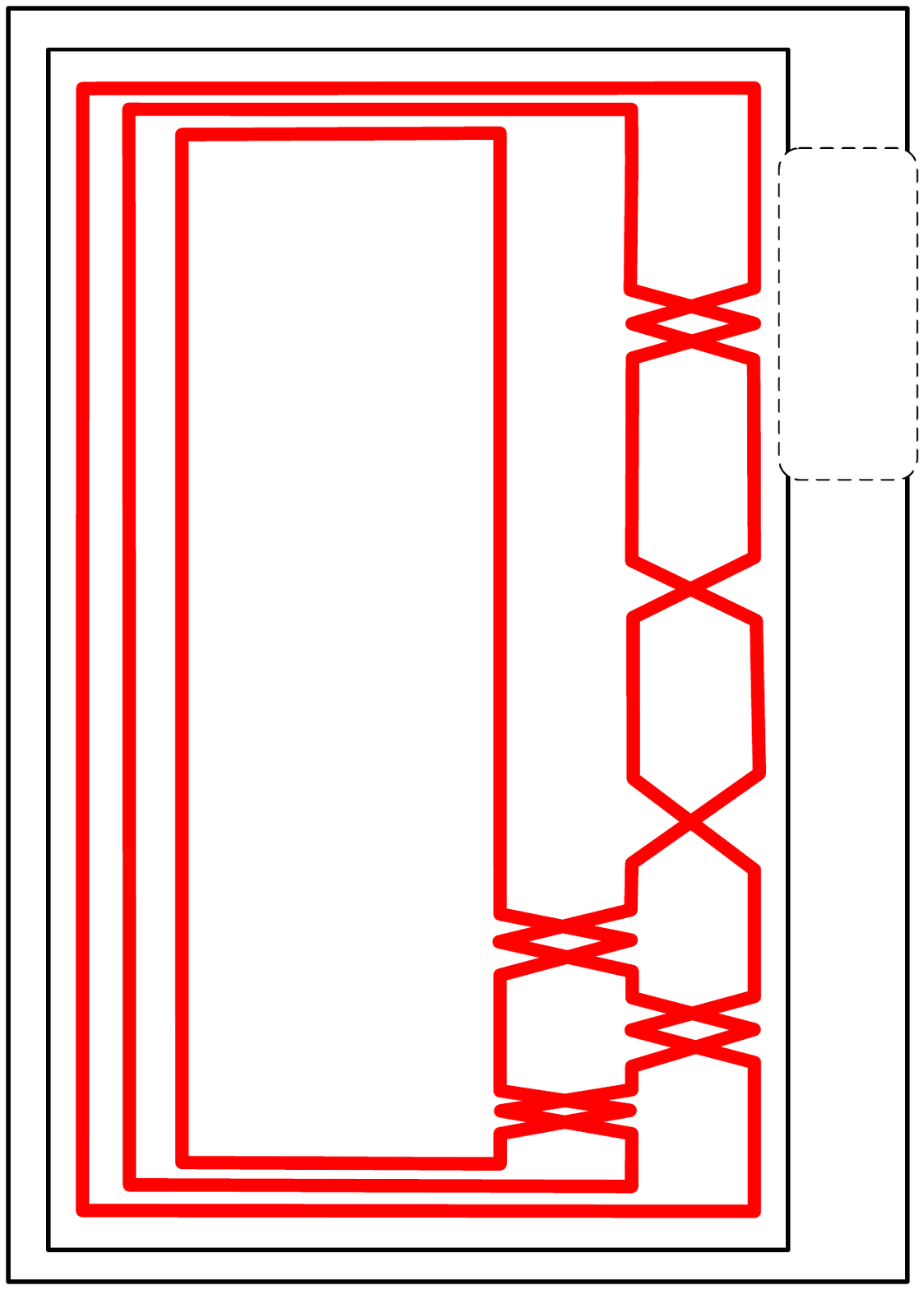}}
            \put(12.2,6.5){$(\text{B})$}
            \put(12.2,4){$\cdots$}
        \end{picture}
        \caption{A stable map $f : S^{3} \to \mathbb{R}^{2}$ with $S_{0}(f)$ isotopic to $K=10_{128}$. The two outermost strands concern the generator $\sigma_1$. The braid word should be read from top to bottom.}
        \label{tegaki8}
    \end{figure}

    \begin{figure}[htbp]
        \setlength\unitlength{1truecm}
        \begin{picture}(15,5)(0,0)
            \put(4.5,0){\includegraphics[width=0.3\textwidth,clip]{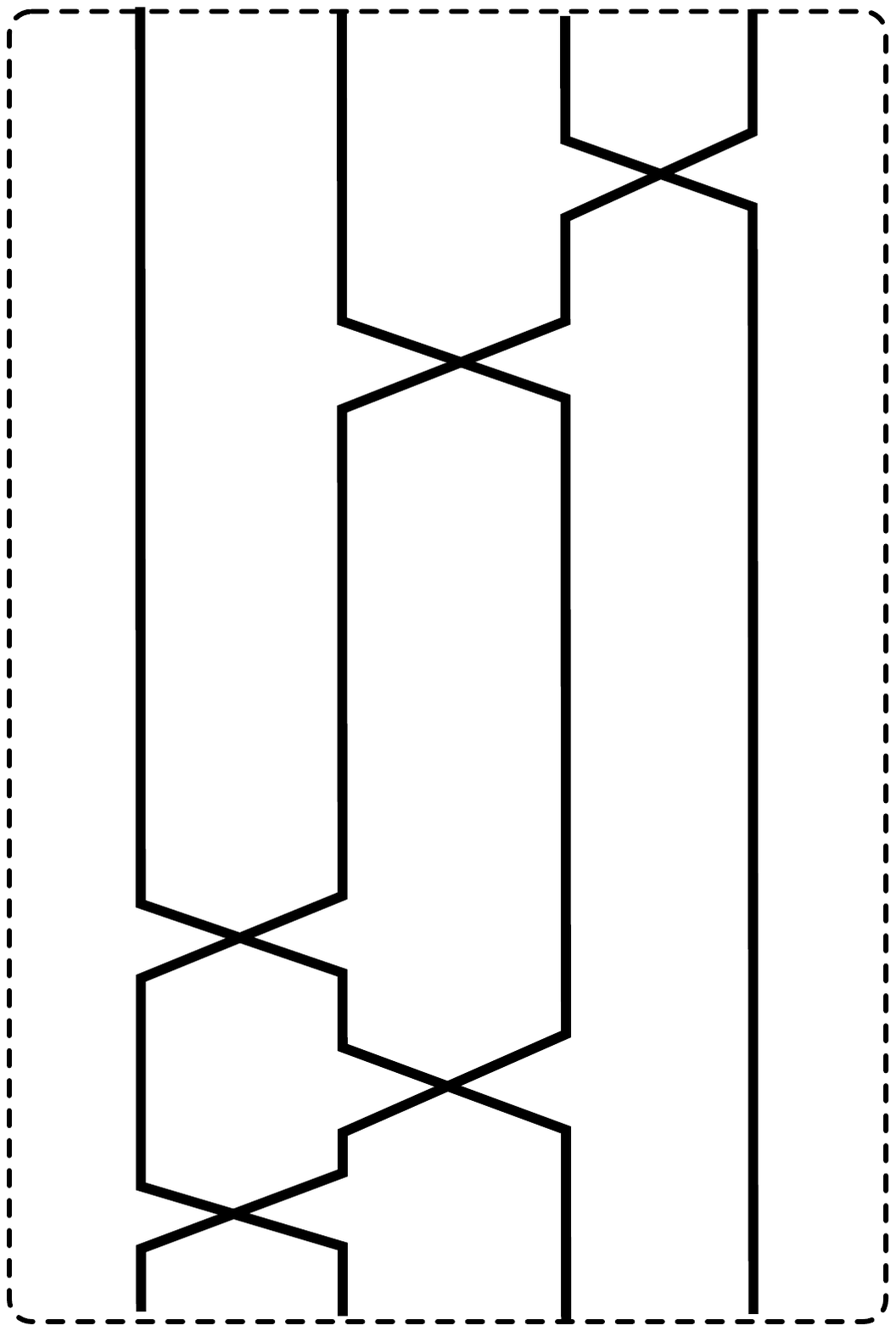}}
        \end{picture}
        \caption{Crossings in a box $(\text{B})$.}
        \label{Box}
    \end{figure}

\end{example}

%%%%%%%%%%%%%%%%%%%%%%%%%%%%%%%%%%%%%%%%%%%%%

\section{a stable map on a $3$-manifold}
\label{s3}

In this section, we give a proof of Theorem~\ref{main2}. Before proceeding to the proof of Theorem~\ref{main2}, we recall the definition of Dehn surgery. Let $K$ be a knot in $S^{3}$, and let $E(K)$ denote the exterior of $K$, defined as $S^{3} \setminus \text{int}\,N(K)$. The isotopy class of an unoriented simple closed curve on $\partial E(K)$ is called a $\textit{slope}$. By using the standard meridian-longitude system, the set of slopes on $\partial E(K)$ is identified with $\mathbb{Q} \cup \{ 1/0 \}$. A $3$-manifold is obtained by gluing a solid torus $V$ to $E(K)$ so that the meridian of $V$ is identified with a curve of slope $p/q$ on $\partial E(K)$. This operation is called $p/q$-$\textit{Dehn surgery}$. Similarly, Dehn surgery on a link is defined for any link. In the following, we consider the case of $p$-Dehn surgery (i.e., $q = 1$), which we call \textit{integral} Dehn surgery. See \cite[Chapter 12]{Lickorish-KT} for details.

\begin{proof}[Proof of Theorem~\ref{main2}.]

It is known that every closed orientable $3$-manifold is obtained from $S^{3}$ by integral Dehn surgery on the closure of a pure braid. See \cite[p.115]{L.S-R.C.}, for example. Suppose that a link $L = K_{1} \cup \cdots \cup K_{l}$ in $S^{3}$ is represented as the closure of a pure braid $b = \sigma_{z_{1}}^{m_{1}} \sigma_{z_{2}}^{m_{2}} \cdots \sigma_{z_{l}}^{m_{l}}$. Since $b$ is a pure braid, the number of components of $L$ is equal to the number of strands of $b$. Then, by Theorem~\ref{main1}, we have a stable map $f : S^{3} \to \mathbb{R}^{2}$ such that $f$ has no cusp points and no singular fibers of type $\mathrm{I\hspace{-1.2pt}I^{3}}$, $S_{0}(f)$ is isotopic to $L$, and $|\mathrm{I\hspace{-1.2pt}I^{2}}(f)| = 2(l-|X|)$ holds with $X = \{ i \mid z_{i}=1,~ 1 \le i \le l \}$. Note that for each $i$, the image $f(K_{i})$ of $K_{i} \subset L$ is a simple closed curve in $\mathbb{R}^{2}$.

Let $D_{1}$ be the image $f(S^{3}) \subset \mathbb{R}^{2}$ and $N(L)$ a sufficiently small closed neighborhood of $L$. By removing the interior of $N(L)$ from $S^{3}$, we obtain the exterior $E(L)$, whose image $f(E(L))$ is a disk $D_{2}$ satisfying $D_{2} \subset D_{1}$. Note that $f(\partial E(L))$ is the union $f(\partial N(K_{1})) \cup \cdots \cup f(\partial N(K_{l}))$, and the image $f(\partial N(K_{i}))$ is a simple closed curve for each $i$. Note that the preimage of a point on $f( \partial N(K_{i}))$ is a meridian of $N(K_{i})$ under $f$. Let $U_{i}$ be a solid torus for each $i \in \{1, \ldots, l\}$. We glue $U_{i}$ to $E(L)$ via integral Dehn surgery on $K_{i} \subset L$ to obtain a $3$-manifold $M$. We remark that a longitude of $U_{i}$ is identified with a longitude of $N(K_{i})$ under the integral Dehn surgery. Then, a smooth map $f_{0} : M \to S^{2}$ is obtained by connecting $f \mid_{E(L)}$ and the natural projections $h_{i} : U_{i} \to \mathbb{D}^{2}$ so that $h_{i}(\partial U_{i}) = f(\partial N(K_{i}))$. This map $f_{0}$ has no cusp points, and its singular set consists of indefinite fold points only. Moreover, $f_{0}$ satisfies the global conditions (5) and (6). Consequently, we obtain a stable map $f_{0} : M \to S^{2}$ such that $f_{0}$ has no cusp points, no definite fold points, and no singular fibers of type $\mathrm{I\hspace{-1.2pt}I^{3}}$, and $|\mathrm{I\hspace{-1.2pt}I^{2}}(f_{0})| = 2(l-|X|)$ holds.

\end{proof}
%%%%%%%%%%%%%%%%%%%%%%%%%%%%%%%%%%%%%%%%%%%%%

\begin{example}

We consider the minimally twisted five-chain link, $L10_{n}113$, in Figure~\ref{tegaki12}. The link is well known in the study of exceptional Dehn surgery. See \cite{M.B-P.C-R.F}, for example. The link $L$ is represented by a $5$-string braid with the braid word $\sigma_{1}^{2} \sigma_{2}^{-2} \sigma_{4}^{-2} \sigma_{3}^{-1} \sigma_{2} \sigma_{1}^{-2} \sigma_{2}^{-1} \sigma_{4}^{2} \sigma_{3}$. Let $M$ be a closed orientable $3$-manifold obtained from $S^{3}$ by integral Dehn surgery on $L$. Then, a stable map $f : M \to S^{2}$ with $S_{0}(f)=\emptyset$ is obtained by the construction given in the proof of Theorem~\ref{main2}. In Figure~\ref{tegaki12}, $f ( S_{1}(f))$ is shown by immersed curves on $S^{2}$. This $f$ has no singular fibers of type $\mathrm{I\hspace{-1.2pt}I^{3}}$. Only some of the singular fibers are shown in the figure.

    \begin{figure}[htbp]
        \setlength\unitlength{1truecm}
        \begin{picture}(15,7)(0,0)
            \put(3,0){\includegraphics[width=0.45\textwidth,clip]{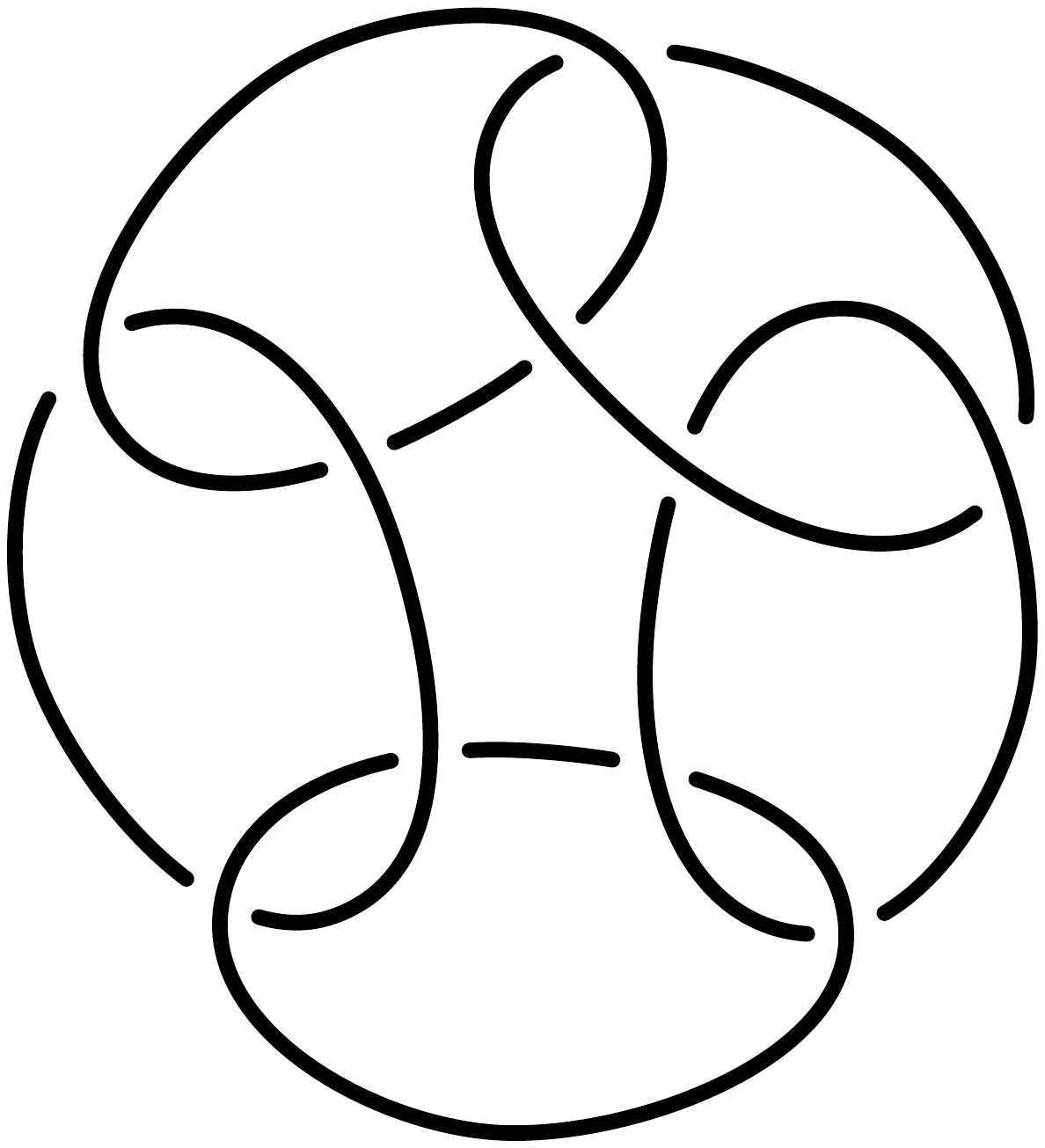}}
        \end{picture}
        \caption{ The link $L10_{n}113$}
        \label{L10ori.}
    \end{figure}

    \begin{figure}[htbp]
        \setlength\unitlength{1truecm}
        \begin{picture}(15,9)(0,0)
            \put(0,0.2){\includegraphics[width=0.45\textwidth,clip]{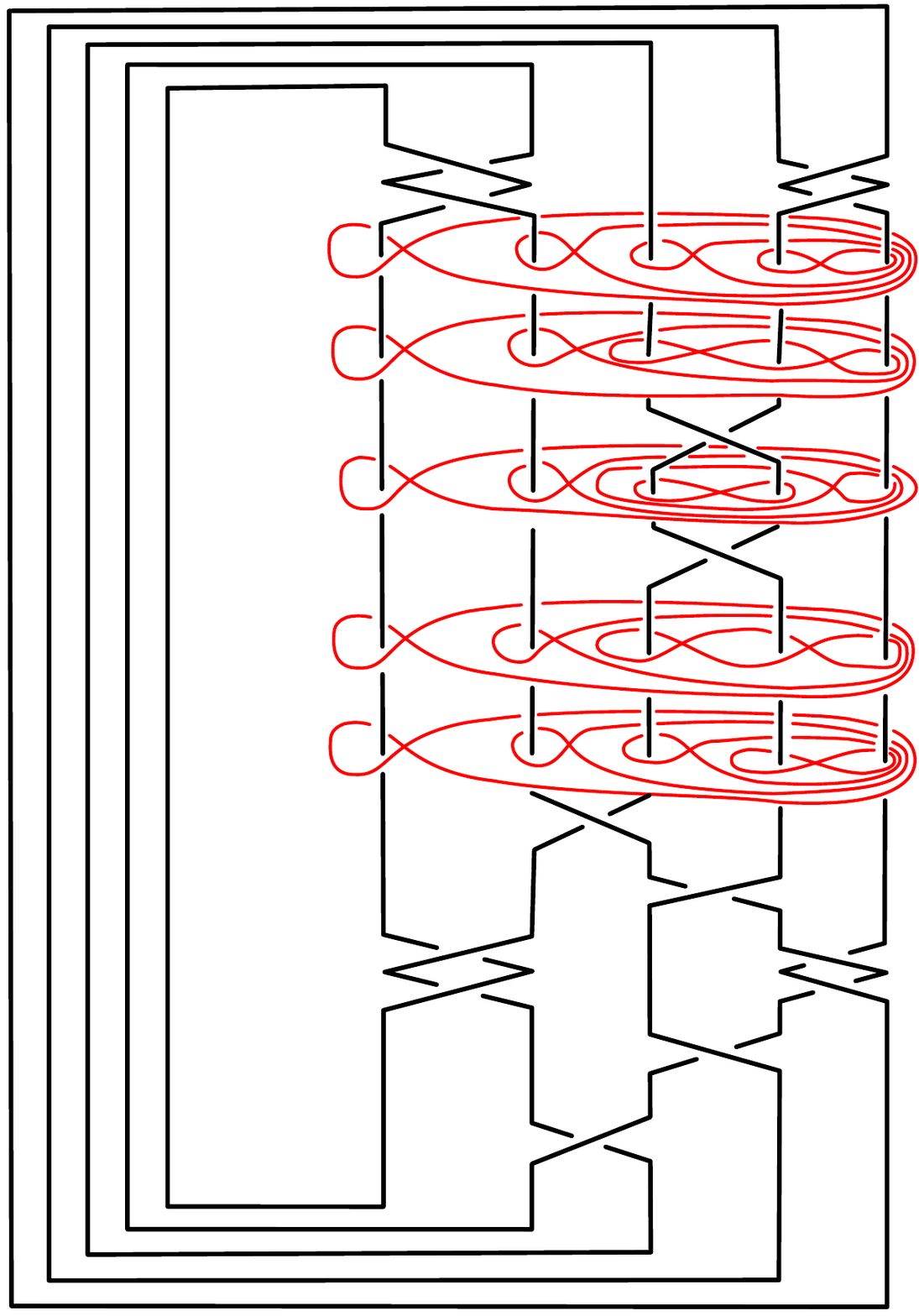}}
            \put(5.8,4.5){$\longrightarrow$}
            \put(6.5,0){\includegraphics[width=0.45\textwidth,clip]{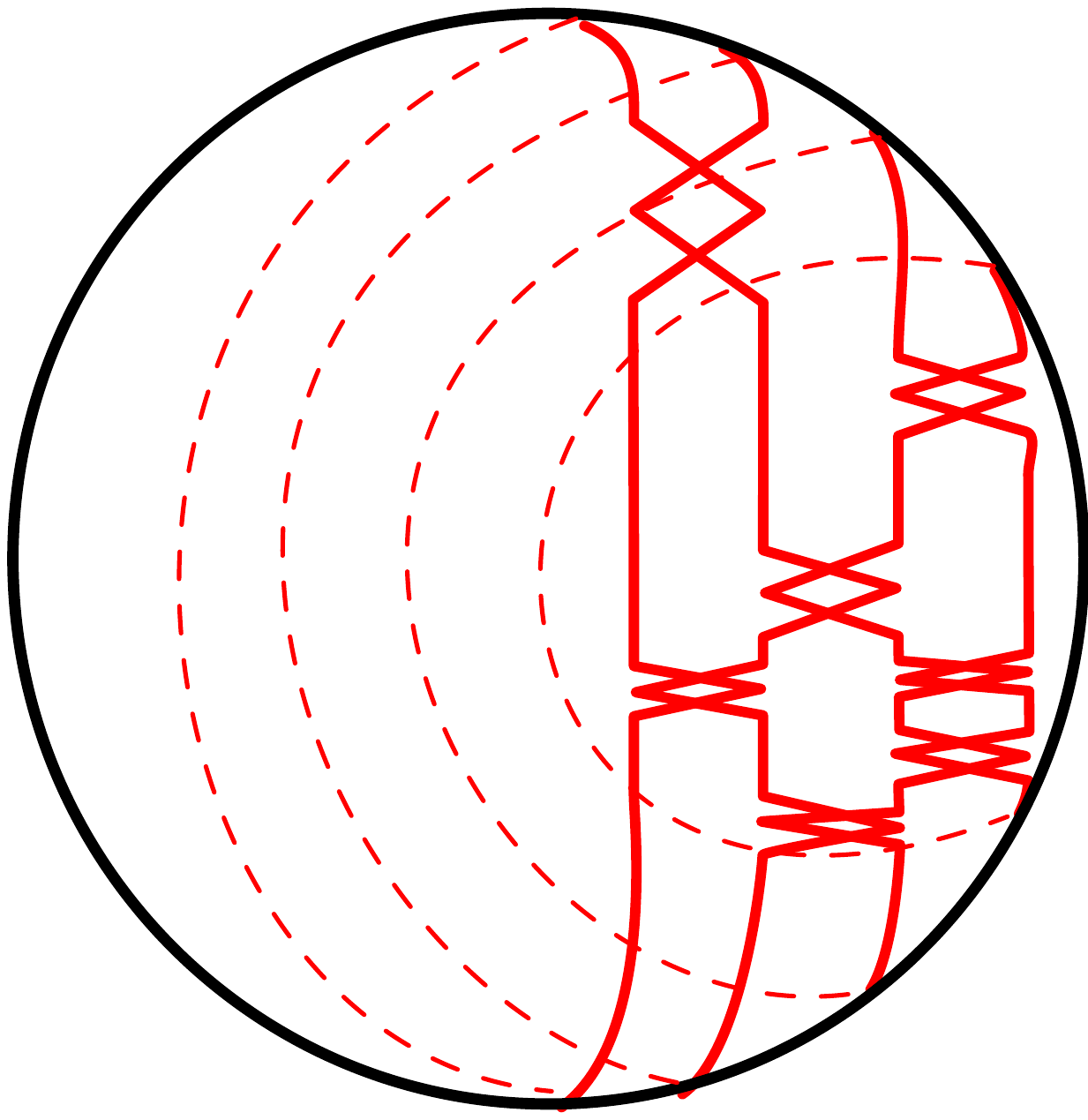}}
        \end{picture}
        \caption{A stable map into the $2$-sphere from a $3$-manifold obtained by integral Dehn surgery on $L10_{n}113$. The generator $\sigma_1$ corresponds to the rightmost twists. The braid word is read from top to bottom.}
        \label{tegaki12}
    \end{figure}
    
\end{example}

%%%%%%%%%%%%%%%%%%%%%%%%%%%%%%%%%%%%%%%%%%%%%

\section*{Acknowledgements}
The author would like to thank his supervisor, Kazuhiro Ichihara, for his continuous guidance, invaluable advice, and warm encouragement.
He also thanks Takahiro Yamamoto for useful comments.

%%%%%%%%%%%%%%%%%%%%%%%%%%%%%%%%%%%%%%%%%%%%%

\bibliographystyle{amsplain}
% \bibliography{name}

\providecommand{\bysame}{\leavevmode\hbox to3em{\hrulefill}\thinspace}
\providecommand{\MR}{\relax\ifhmode\unskip\space\fi MR }
\providecommand{\MRhref}[2]{%
  \href{http://www.ams.org/mathscinet-getitem?mr=#1}{#2}
}
\providecommand{\href}[2]{#2}

\end{document}